\documentclass[10.5pt,a4paper]{article}

\usepackage{amssymb}
\usepackage{amsmath}
\usepackage{tensor}
\usepackage{color}
\usepackage{stmaryrd}
 \usepackage{relsize}
\usepackage{cite}
\usepackage[utf8]{inputenc}
\newtheorem{Th}{Theorem}[section]
\newtheorem{Prop}{Proposition}[section]
\newtheorem{Lm}{Lemma}[section]

\newtheorem{Rm}{Remark}

\newcommand{\be}{\begin{equation}}
\newcommand{\ee}{\end{equation}}
\newcommand{\bes}{\begin{equation*}}
\newcommand{\ees}{\end{equation*}}
\newcommand{\R}{\mathbb{R}}

\newcommand\res{\mathop{\hbox{\vrule height 7pt width .5pt depth 0pt
\vrule height .5pt width 6pt depth 0pt}}\nolimits}

\def\theequation{\thesection.\arabic{equation}}
\def\theTh{\Roman{section}.\arabic{Th}}
\def\theProp{\Roman{section}.\arabic{Prop}}
\def\theCo{\Roman{section}.\arabic{Co}}
\def\theLm{\Roman{section}.\arabic{Lm}}

\def\theRm{\Roman{section}.\arabic{Rm}}
\newcommand{\reset}{\setcounter{equation}{0}\setcounter{Th}{0}\setcounter{Prop}{0}\setcounter{Co}{0}\setcounter{Lm}{0}\setcounter{Rm}{0}}

\def\al{\alpha}
\def\la{\lambda}
\def\eps{\varepsilon}

\def\Om{B}
\def\om{\omega}

\def\pro{\pi_{\vec{n}}}

\def\bh{\vec{h}}

\def\bH{\vec{H}}
\def\bC{\vec{C}}

\def\bA{\vec{A}}
\def\bF{\vec{F}}
\def\bU{\vec{U}}

\def\bL{\vec{L}}

\def\bR{\vec{R}}

\def\bV{\vec{V}}

\def\bq{\vec{q}}

\def\bv{\vec{v}}

\def\bp{\vec{\Phi}}

\def\bT{\vec{T}}

\def\di{{D}}

\def\res{\mathop{\hbox{\vrule height 7pt width .5pt 
depth 0pt\vrule height .5pt width 6pt depth 0pt}}\nolimits}

\begin{document}

\title{Analysis of Critical Points of\\ Conformally Invariant Curvature Energies in 4d}
\author{Yann Bernard\footnote{School of Mathematics, Monash University, 3800 Victoria, Australia.}}
\date{\today}
\maketitle

{\bf Abstract:} {\it We consider critical points of a conformally invariant energy of the type
$$
\int_{\Sigma}|\pro d\bH|^2-|\bH\cdot\bh|^2+7|\bH|^4+\mathcal{O}(|\bh_0|^4)\:,
$$  
for a submanifold $\Sigma^4\subset\R^{m\ge5}$. We show that assuming $\bp\in W^{3,2}\cap W^{1,\infty}$ is sufficient to give a meaning to the Euler-Lagrange equation $\vec{\mathcal{W}}=\vec{0}$, which is a quasilinear system of sixth-order equations in the immersion $\bp$. Using a reformulation ultimately stemming from Noether's theorem, we show that weak solutions are smooth and we deduce an $\eps$-regularity estimate. We also consider local Palais-Smale sequences and show they converge to a solution of the constrained equation $\vec{\mathcal{W}}=-\langle T,\bh\rangle$, for some divergence-free traceless symmetric 2-tensor $T$ acting as a Lagrange multiplier. }\\



\reset

\section{Introduction and Main Results}


Conformally invariant equations arise naturally in differential geometry. These are equations, often of elliptic type, whose form is preserved under conformal rescaling of the metric \( g \mapsto e^{2\omega}g \). Such invariance reflects geometric structure and provides connections between curvature, topology, and analysis. The study of these problems has been central to understanding scalar and higher-order curvature invariants, moduli of conformal structures, and critical phenomena in both mathematics and physics. A classical example is the Yamabe problem, which asks whether a compact Riemannian manifold \( (M^n, g) \), \( n \geq 3 \), admits a conformal metric with constant scalar curvature. This reduces to solving the nonlinear PDE, which is conformally invariant in the sense that both the operator and the scalar curvature transform appropriately under conformal changes. The Yamabe problem was resolved through the combined efforts of Yamabe, Trudinger, Aubin, and Schoen \cite{yamabe1960trials,trudinger1968remarks,aubin1976equations,schoen1984conformal}. It has served as a model for many subsequent conformally invariant variational problems.\\
In dimension four, the role of the scalar curvature is taken up by the so-called Q-curvature, a fourth-order scalar invariant. The associated Paneitz operator \( P_g \), which is conformally covariant under changes \( g \mapsto e^{2\omega}g \), generalizes the conformal Laplacian. The Paneitz equation is a natural analogue of the Yamabe equation in four dimensions and features prominently in the study of conformal anomalies and extremal metrics in conformal classes \cite{paneitz2008quartic,branson1995sharp,chang1995prescribing}. The relevance of Q-curvature emerges from the volume renormalisation of conformally compact manifolds, see \cite{FG, HS, W}. 
Further in dimension four, the Bach tensor \( B_{ij} \) arises as the Euler--Lagrange equation of the conformally invariant functional
\[
\mathcal{E}_B(g) = \int_M |W_g|^2 \, d\mu_g,
\]
where \( W_g \) is the Weyl tensor. The Bach tensor is given by
\[
B_{ij} = \nabla^k \nabla^\ell W_{ikj\ell} + \tfrac{1}{2} R^{k\ell} W_{ikj\ell},
\]
and satisfies \( B = 0 \) if and only if \( g \) is a critical point of \( \mathcal{E}_B \). It is symmetric, trace-free, and conformally invariant. All Einstein metrics in 4d are Bach-flat. The Bach equation plays a central role in conformal geometry and mathematical relativity.\\
Fully nonlinear conformal equations have also attracted considerable attention. A key example is the \( \sigma_k \)-Yamabe problem, which asks for conformal metrics with prescribed \( \sigma_k \)-curvature -- the \( k \)-th elementary symmetric function of the eigenvalues of the Schouten tensor. These equations are generally fully nonlinear and degenerate elliptic, but they retain conformal invariance under suitable geometric conditions. Their study links conformal geometry with nonlinear potential theory and has produced many deep results \cite{viaclovsky2000conformal,guan2003local}.\\

Other notable examples of conformally invariant problems arise from submanifold geometry, particularly in the study of surface immersions in Euclidean or spherical space. One of the most significant conformally invariant variational problems in this setting is the Willmore energy, defined for a surface \( \Sigma^2 \subset \mathbb{R}^{m\ge3} \) by
\[
\mathcal{W}(\Sigma) = \int_\Sigma |\bH|^2 \, d\mu,
\]
where \( \bH \) is the mean curvature vector of \( \Sigma \) and \( d\mu \) is the induced area form. This energy is invariant under Möbius transformations, and its critical points include minimal surfaces and round spheres. The Willmore conjecture, proposed in the 1960s, stated that among all tori immersed in \( \mathbb{R}^3 \), the Clifford torus minimizes \( \mathcal{W} \), achieving the value \( 2\pi^2 \). This was resolved affirmatively by Marques and Neves using min-max theory \cite{MN}. The critical points of the Willmore energy -- be they minimisers or not -- are of course interesting in their own right. They are smooth immersions $\bp:\Sigma\rightarrow\mathbb{R}^m$ satisfying a quasilinear elliptic equation of order four (for the immersion):
\be\label{will2}
\Delta_\perp\bH+(\bH\cdot\bh^{ij})\bh_{ij}-2|\bH|^2\bH\;=\;\vec{0}\:.
\ee
Here $\bh$ is the second fundamental form, $\Delta_\perp:=\pro\nabla_j\pro\nabla^j$ is the covariant Laplacian in the normal bundle, and $\pro$ denotes projection onto the normal space. Repeated indices are summed over using the induced metric $g_{ij}=\nabla_i\bp\cdot\nabla_j\bp$. We note that minimal surfaces and their conformal transforms are examples.\\

\noindent
The analysis of solutions to (\ref{will2}) is particularly challenging. The fourth-order operator, of course, precludes the use of maximum principle. Moreover, under the natural hypothesis that $\bh$ be square-integrable is largely insufficient to give a weak meaning to the nonlinearities of the Willmore equation. To thwart these difficulties in codimension 1, Kuwert and Sch\"atzle devised a powerful ``ambient" approach based on repeated integration by parts and estimates from the Michael-Simon-Sobolev inequality \cite{KS1, KS2}. A different approach was independently proposed by Tristan Rivi\`ere within the framework of weak immersions \cite{Riv1}. In this ``parametric" method, codimension is irrelevant, and one starts by observing that (\ref{will2}) is in divergence form -- a fact which ultimately follows from Noether's theorem and the translation invariance of the Willmore energy \cite{Ber1}. Two more invariant equations are found from the other invariances (rotation and dilation). Integrating once these three conservation laws yields three potentials. Rivi\`ere's tour-de-force consists in relating these potentials in a system of PDEs with good analytical dispositions, and in particular permitting the use of integration by compensation techniques. A brief summary of Rivi\`ere's formulation of the Willmore problem is given in the introduction of \cite{Ber2}. His original framework has inspired numerous subsequent works by various authors. We content ourselves with citing \cite{BR1, BR2, BLM, LR, PR, Riv2}, and caution the reader that many relevant references are omitted from this deliberately short list. \\

Conformally invariant problems are significant not only for their geometric content but also for their analytical challenges. These include the absence of a maximum principle for higher-order operators, the appearance of critical Sobolev exponents leading to loss of compactness, and the presence of large symmetry groups that complicate uniqueness and stability analysis. Blow-up phenomena, bubbling, and concentration of curvature are typical features, necessitating intricate compactness and rescaling techniques. Applications range from classification results and rigidity theorems in differential geometry to physical theories such as conformal field theory and the AdS/CFT correspondence, where renormalised conformal invariants play a key role. Broadly speaking, the {\it Anti-de-Sitter/Conformal Field Theory correspondence} of Maldacena \cite{Mal} postulates the existence of a duality equivalence between gravitational theories (such as string or $M$ theory) on anti-de Sitter spaces $M$ and conformal field theories on the boundary at conformal infinity $\partial M$. This correspondence and the search for a rigorous mathematical formulation has led to some surprising discoveries in geometry, and in particular yields new purely geometric invariants of conformally compact Einstein manifolds and of their minimal submanifolds. The relevant notion of conformal infinity is that introduced by Penrose. A Riemannian metric $g_+$ on the interior $X^{n+1}$ of a compact manifold with boundary $\bar{X}$ is said to be conformally compact if $\bar{g}=r^2g_+$ extends continuously (or with some degree of smoothness) as a metric to $\bar{X}$, where $r$ is a defining function for $M=\partial X$, i.e. $r>0$ on $X$ and $r= 0$, $dr= 0$ on $M$. The restriction of $\bar{g}$ to $TM$ rescales upon changing $r$, so defines a conformal class of metrics on $M$, the conformal infinity of $g_+$. Of particular interest  are conformally compact metrics $g_+$ which are Einstein. 
According to the AdS/CFT correspondence, in a suitable approximation, the expectation value of an observable for a submanifold $N$ of $M$ can be calculated in terms of the area $A(Y )$ in the $g_+$ metric of minimal submanifolds $Y$ of $X$ with $\partial Y = N$. Necessarily $A(Y ) = \infty$, so one is led to consideration of renormalizing the area of a minimal submanifold. If $r$ is the special defining function associated to a conformal representative on $M$, then $\text{Area}(Y \cap \{r > \eps\})$ has an expansion in negative powers of $\eps$, and a $\log\eps$ term if $k = \text{dim}(N)$ is even. The $\log\eps$ coefficient is a conformal invariant of the submanifold $N$ of $M$. For 2-dimensional $N$, this coefficient is a version on a general conformal manifold of the Willmore functional of a surface in conformally flat space, of which minimal submanifolds are critical points. The Willmore functional is called the “rigid string action” in the physics literature. The same procedure was independently carried out by Graham and Reichert in \cite{GR} and by Zhang in \cite{Zha} for four-dimensional submanifolds of $\mathbb{R}^{m\ge5}$. The authors found the energy
$$
\mathcal{E}_A\;:=\;\int_\Sigma\big(|\pro d\bH|^2-|\bH\cdot\bh|^2+7|\bH|^4\big) d\text{vol}_g
$$
turns out to be conformally invariant. In codimension 1, this energy was first identified by Guven \cite{Guv}; which the author derives from scratch by evaluating and cleverly recasting the defect of conformal invariance of the Dirichlet energy of mean curvature. Blitz, Gover, and Waldron \cite{BGW} use instead tractor calculus and produce numerous conformal invariant energies, including all those which appear in the present paper. \\

The Euler-Lagrange equation for $\mathcal{E}_A$ is derived in \cite{Ber2} (see also \cite{GR} and \cite{WY}):
\begin{eqnarray}\label{horrib}
&&\vec{\mathcal{W}}\;\;=\;\;\dfrac{1}{2}\Delta_\perp^2\bH+\dfrac{1}{2}\Delta_\perp\langle\bH\cdot\bh,\bh\rangle-7\Delta_\perp(|\bH|^2\bH)+8\pro\nabla_j\big(\bH\nabla^j|\bH|^2\big)+4\pro\nabla_j\big((\bH\cdot\bh^j_i)\pro\nabla^i\bH\big)\nonumber \\
&&\hspace{1.5cm}+\:\Big[\dfrac{1}{2}\bh^{ij}\cdot\Delta_\perp\bH-2\nabla_i\bH\cdot\pro\nabla_j\bH+2(\bH\cdot\bh^{ik})(\bH\cdot\bh^j_k)  \nonumber\\
&&\hspace{2cm}+\:\dfrac{1}{2}\langle\bH\cdot\bh,\bh\rangle\cdot \bh^{ij}-7|\bH|^2\bH\cdot\bh^{ij} +\big(|\pro d\bH|^2-|\bH\cdot\bh|^2+7|\bH|^4\big)g^{ij}\Big]\bh_{ij}\nonumber\\
&&\hspace{.5cm}\;\;=\;\;\vec{0}\:.
\end{eqnarray}
The energy $\mathcal{E}_A$ is, for many of its aspects, the analogue of the two-dimensional Willmore energy. Both are conformally invariant, both include minimal immersions as critical points, both have a quasilinear elliptic system (of order 4 in the immersion for 2d Willmore, and of order 6 in the four-dimensional setting). However, unlike the Willmore energy, $\mathcal{E}_A$ is not bounded from below. Recent works have shown that minimisers of $\mathcal{E}_A$ are hard to come by. Martino \cite{Mart} proves that $\mathcal{E}_A$ is unbounded from below on closed hypersurfaces. In an effort to thwart this annoying feature, one may add some other lower-order conformally invariant energy. For example, as is elementary to check \cite{Vya}, the energy $\mathcal{E}_A+\beta\int|\bh_0|^4$ is bounded below by $8\pi^2$ for $\beta>\frac{1}{12}$. In the recent work \cite{WY}, the authors establish a connected sum energy reduction for the energy $\mathcal{E}_A$. In the submanifold setting, the diversity of conformally invariant energies increases, if only due to combinatorics. Indeed, from the pointwise conformal invariant tensor $\bh_0$ one can construct four distinct conformally invariant energies with integrands 
$$
|\bh_0|^4\quad,\quad \langle\bh_0^4\rangle:=\big(\bh_0^{ij}\cdot\bh_0^{kl}\big)\big((\bh_0)_{ij}\cdot(\bh_0)_{kl}\big)
$$
and
$$
|h_0^2|^2:=\big((\bh_0)^{ik}\cdot(\bh_0)_k^j\big)\big((\bh_0)_{il}\cdot(\bh_0)_j^l\big)\quad,\quad \text{Tr}(\bh_0^4):=\big((\bh_0)^{ij}\cdot(\bh_0)^{kl}\big)\big((\bh_0)_{il}\cdot(\bh_0)_{kj}\big)\:.
$$
In codimension 1, the four energies coincide pairwise. \\

We will be concerned with energies of the type $\mathcal{E}=\mathcal{E}_A+\mathcal{E}_0$, where $\mathcal{E}_0$ denotes any linear combinations of the four conformally invariant energies given above. In \cite{Ber2}, we show how to use Noether's theorem to develop an analytically favourable reformulation of the Euler-Lagrange equation. This reformulation will enable us to prove the regularity of weak solutions as well as to consider the limiting behaviour of a Palais-Smale sequence. \\

We note that, while critical points of $\mathcal{E}_A+\mathcal{E}_0$ satisfy a quasilinear sixth-order system of equations, critical points of $\mathcal{E}_0$ satisfy a fully nonlinear system of equations. The methods developed in the present work are unsuitable to study critical points of energies of the type $\mathcal{E}_0$. One prominent linear combination of the type $\mathcal{E}_0$ is the intrinsic {\it Bach energy} appearing in the introduction \cite{Ber2}. Its (intrinsic) critical points are known as {\it Bach flat} manifolds and include Einstein metrics as examples. Their analytical properties are discussed in \cite{TV}. \\

The present paper presents an analytically suitable framework to obtain critical energy estimates. Much like Rivi\`ere's parametric approach initiated numerous works on the two-dimensional Willmore problem (see introduction and references therein), we anticipate the setting given in the present paper will yield analogous results for four-dimensional conformally invariant curvature energies, and in time produce a better picture of their energy landscapes. In particular, since conformal transforms of minimal submanifolds are critical points of $\mathcal{E}_A$, it might be possible to study complete, non-compact minimal immersions, much like what was done in \cite{BR3}. This will require a solid understanding of the behaviour of critical points of $\mathcal{E}_A$ near a singular point.

\subsection{Regularity of critical points}\label{introreg}


We will be concerned with non-degenerate immersions $\bp:\Omega\rightarrow\R^{m\ge5}$. We ask that the pull-back metric $g$ be uniformly elliptic. In order to give a meaning to the energy $\mathcal{E}$, we also ask that $\bp\in W^{3,2}\cap W^{1,\infty}$. This guarantees the metric $g_{ij}$ lies in $W^{2,2}\cap L^\infty$.\\

\noindent
In this context, it appears difficult to give a meaning to the Euler-Lagrange equation (\ref{horrib}). Its right-hand side does not even lie in $L^1$, as it contains, for example, terms of order $\mathcal{O}(|\bh|^5)$. A similar phenomenon appears in the two-dimensional Willmore equation $\Delta_\perp\bH+2(\bH\cdot\bh_0)\bH=\vec{0}$. Imposing $\bh\in L^2$ does not give any meaning to the left-hand side. \\

\noindent
Of course, from the point of view of {\it a priori} estimates, one may suppose first that $\bp$ is smooth. Via repeated multiplication by astutely chosen test-functions and integration by parts, it is possible to derive ``ambient" energy estimates, as done in \cite{KS1, KS2} (only in codimension 1). In dimension four, the same approach might be possible, but would unmistakably be particularly tedious, due to the high number of derivatives involved and the variety of nonlinear terms. One would most likely require improved Gagliardo-Nirenberg inequalities for hypersurfaces, such as those given in \cite{FHT}. \\

\noindent
We favour in this work a different approach inspired by Rivi\`ere's study of Willmore immersions in dimension 2 \cite{Riv1}. This approach consists first in replacing the Euler-Lagrange equation $\vec{\mathcal{W}}=\vec{0}$ by $d^{\star_g}\bV=\vec{0}$, for a suitable vector-field $\bV$, which can be determined \cite{Ber2}. This reformulation of the Euler-Lagrange equation has a weak sense for $\bp\in W^{3,2}\cap W^{1,\infty}$ and is in fact sufficient to initiate a regularity process. Our proof applies regardless of codimension. \\

 As the metric is assumed to be non-degenerate, we deduce that $|g|^{1/2}$, $g^{ij}$ belong to $W^{2,2}\cap L^\infty$. Since $\bH=g^{ij}\partial_{ij}\bp$ and since $W^{2,2}\cap W^{1,\infty}$ preserves continuously the space $W^{1,2}$ by multiplication, we find that $\bH$ belongs to $W^{1,2}$. Similarly, the second fundamental form $\bh$ lies in $W^{1,2}$. Again, the position of the indices in $\bh$ is irrelevant, since the metric multiplies $W^{1,2}$. By the Sobolev-Lorentz embedding theorem (see Appendix), we obtain that $\bh$ lies in the Lorentz space $L^{4,2}$. These considerations bring us to the following smallness hypothesis
\be\label{epsreg0}
\Vert\bh\Vert_{L^{4,2}(\Omega)}\;<\;\eps_0\:,
\ee
for some $\eps_0>0$ which will be adjusted as necessary. By simple scaling, we can without loss of generality suppose that (\ref{epsreg0}) holds. 

\begin{Th}\label{Threg}
Let $\Omega$ be a ball in $\mathbb{R}^4$, and let $\bp:\Omega\rightarrow\R^{m\ge5}$ be a non-degenerate immersion in $W^{3,2}\cap W^{1,\infty}(\Omega)$ with pull-back metric $g=\bp^\star g_{\R^{m\ge5}}$. We suppose that $g$ is non-degenerate: there exists a constant $c_0>0$ such that
$$
c_0\,\delta_{ij}\;\leq\; g_{ij}\;\leq\;\dfrac{1}{c_0}\delta_{ij}\:.
$$
Let $\bp$ be a critical point of the energy $\mathcal{E}$. Then:\\[1ex]
There exists $\eps_0>0$ such that if (\ref{epsreg0}) holds, then $\bp$ is smooth in the interior of $\Omega$. Moreover, the following estimate holds
$$
\dfrac{1}{r}\Vert\bh\Vert_{L^{\infty}(B_r)}+\Vert \pro d\bH\Vert_{L^\infty(B_r)}\;\lesssim\;\dfrac{1}{r^2}\Big(\Vert\bh\Vert_{L^4(B)}+\Vert \pro d\bH\Vert_{L^2(B)}\Big)\:,
$$
up to a multiplicative constant depending on $c_0$, $\eps_0$, and $\Vert g\Vert_{W^{2,2}\cap L^\infty(\Omega)}$. In this estimate, $B$ is a compact ball properly contained in $\Omega$, while $B_r$ denotes $B$ rescaled by a factor $r\in(0,1)$.
\end{Th}
$\hfill\square$\\

\noindent
An important ingredient in the proof of Theorem~\ref{Threg} is the possibility to Hodge decompose $L^p$ differential forms {\it with natural estimates}. While this theory is well-established when the metric is sufficiently regular (typically $C^1$, \cite{Iwa}), things are far less obvious when the metric is merely bounded. Fortunately, as shown in \cite{Mcintosh}, for certain range of $p$ of the type $p_H\leq p\leq p^H$ with $p_H<2<p^H$, Hodge decompositions with estimates are at reach. The exact values of $p_H$ and $p^H$ will be irrelevant to us, as we will always be concerned with values of $p$ near $2$. Similarly, $L^p$-theory for the Hodge Laplacian $dd^{\star_g}+d^{\star_g}d$ is available for $p$ in the same range, and that will be decisive for us. Finally, we note that the pull-back metric $g$ lies in $W^{2,2}\subset\text{VMO}$. Estimates for solutions of second-order elliptic equations with VMO coefficients are well understood \cite{Chiarenza} and they will also be of help. The analytical superiority of VMO over $L^\infty$ is that $C^{\infty}_c$ is dense in VMO (for the BMO norm). \\

The proof of Theorem~\ref{Threg} relies on a delicate argument involving particular structures with good analytical dispositions. In broad terms, as the energy $\mathcal{E}$ is invariant under translation, and one can call upon Noether's theorem to write the Euler-Lagrange equation in divergence form:
$$
d^\star\bV\;=\;\vec{0}\:.
$$
The exact formulation and properties of the vector-field $\bV$ were studied at length in \cite{Ber2}. It can be computed directly from the immersion $\bp$. Some of these properties will be recalled in the course of our arguments. Integrating once the divergence form yields the existence of a closed form 1-vector-valued 2-form $\bL$ satisfying $d^\star\bL=\bV$. Similarly, the invariance under dilation and rotation give further a scalar-valued 2-form $S$ and a 2-vector-valued 2-form $\bR$. One verifies (for the notational conventions used, please refer to Section~\ref{nota} in the Appendix) that
$$
d^\star S=\bL\stackrel{\!\!\!\cdot}{\res}d\bp\quad,\quad dS=2\bL\stackrel{\cdot}{\wedge}d\bp\quad,\quad d^{\star}\bR=\bL\stackrel{\!\!\!\wedge}{\res}d\bp+\bU\quad,\quad dR=2\bL\stackrel{\wedge}{\wedge}d\bp\:,
$$
with
\be\label{defU}
\bU\;:=\;\dfrac{1}{2}\Delta_\perp\bH\wedge d\bp+\mathcal{O}(|\bh|^3+|\bh||\pro d\bH|\big)\:.
\ee
What is remarkable is that the potentials $S$ and $\bR$ may be recombined so as to eliminate both $\bL$ and $\Delta_\perp\bH$. More precisely, we will show that the following equation holds:
\bes
d^\star\big(3\bR+\vec{\eta}\stackrel{\bullet}{\odot}\bR+\vec{\eta}\odot S\big)+d\big(\langle \vec{\eta}\stackrel{\bullet}{,}\bR\rangle+\langle \vec{\eta},S\rangle\big)\;=\;\mathcal{O}\big(|\bh|^3+|\bh||\pro d\bH|+|\bh||\bR|+|\bh||S|\big)\:.
\ees
In this expression, $\vec{\eta}$ is the 2-vector-valued 2-form $\vec{\eta}=\dfrac{1}{2}d\bp\stackrel{\wedge}{\wedge}d\bp$. Applying $d^\star$ to both sides will enable us to obtain estimates for the 2-vector-valued scalar
$$
\bq\;:=\;\langle \vec{\eta}\stackrel{\bullet}{,}\bR\rangle+\langle \vec{\eta},S\rangle\:.
$$
In turn, this information is used to recover estimates for geometric quantities. We will see that 
\bes
2\Delta_\perp\bH\;=\;\pro(\nabla^i\bq\,\res \nabla_i\bp)+\mathcal{O}\big(|\bh|^3+|\bh||\pro d\bH|+|\bh||\bR|+|\bh||S|\big)\:,
\ees
which warrants a return to geometric from the potentials $S$ and $\bR$. The rest of the regularity proof proceeds more or less as it did in the two-dimensional Willmore situation \cite{Riv1}, first by upgrading $\bH\in W^{1,2}$ to $\bH\in W^{1,s}$, for some $s>2$.\\

Translating this information back to the immersion is done via the equation $\Delta_g\bp=4\bH$. Differentiating gives $\Delta_{\text{Hodge}}d\bp=4d\bH$. Choosing $2<s<p^H$ (with $p^H$ as in \cite{Mcintosh}) guarantees that $d\bp$ lies in $W^{2,s}\subset C^{0,\al}$, for some $\al>0$. The metric regularises accordingly and one may invoke standard elliptic estimates to initiate a bootstrap.

\subsection{Palais-Smale sequences}

We have seen that critical points of our energy satisfy 
$$
d^{\star_g}\bV\;=\;\vec{0}\:,
$$
for some vector-field $\bV$ which is explicitly computed from the immersion $\bp$. For a non-degenerate) metric, this equation is equivalent to 
$$
\partial_i\big(|g|^{1/2}g^{ij}\bV_j\big)\;=\;\vec{0}\:.
$$
In this section we consider a sequence of immersions $\bp_n\in W^{3,2}\cap W^{1,\infty}$ such that the corresponding sequence of induced metrics $g_n$ remains within the same conformal class: $g_n=e^{2\la_n}g_0$, with $g_0$ a smooth reference metric. We demand that $(\la_n)_n$ be uniformly bounded in $W^{2,(2,1)}$. \\

We will be concerned with (non-degenerating) sequences of immersions satisfying the following {\it local Palais-Smale}-type condition\footnote{Prime indicates the dual space.}:
\bes
\partial_i\big(|g_n|^{1/2}g_n^{ij}(\bV_n)_j\big)\:\longrightarrow\:\vec{0}\qquad\text{strongly in}\:\:\:(W^{3,2}\cap W^{1,\infty})'\:,
\ees
The vector-field $\bV_n$ is that associated to $\bp_n$. 

\begin{Th}\label{ThPS}
Let $\bp_n:B\rightarrow\R^{m\ge5}$ be a non-degenerating sequence of immersions uniformly bounded in $W^{3,2}\cap W^{1,\infty}(\Om)$. We ask that the induced metrics $g_n$ all lie in the conformal class of the regular metric $g_0$. We suppose further that the sequence of conformal parameters $\la_n$ be uniformly bounded in the space $W^{2,(2,1)}$. Then:\\

There exist an immersion $\bp\in W^{3,2}\cap W^{1,\infty}(B)$ in the conformal class of $g_0$, and a symmetric trace-free, divergence-free 2-tensor $T$ such that, up to extraction of a subsequence,
$$
\bp_n\;\stackrel{w}{\longrightarrow}\;\bp\qquad\text{in}\:\:W^{3,2}(B)\:,
$$ 
and $\bp$ is a ``constrained" critical point of $\mathcal{E}$ satisfying the perturbed Euler-Lagrange equation
$$
\vec{\mathcal{W}}\;=\;-\langle T,\bh\rangle\qquad\text{in}\:\:\:\mathcal{D}'(\Om)\:.
$$
\end{Th}

Assuming that the sequence of pull-back metrics $g_n$ lies in the same conformal class of a regular metric $g_0$ shifts all regularity questions related to the metric $g_n$ reduce to the regularity of the conformal factor $\la_n$. 
We remark that $W^{2,(2,1)}$ embeds continuously into $W^{2,2}\cap L^{\infty}$. It would seem more natural to assume that $\la_n$ is uniformly bounded in $W^{2,2}\cap L^\infty$. As we will see in the proof, though, it is necessary to suppose the sequence of conformal parameters is uniformly bounded in a slightly smaller space. I was unable to weaken this hypothesis. In dimension 2, matters are quite different since all metrics are conformally flat. The Liouville equation reads $\Delta_{\text{flat}}\la_n=e^{2\la_n}K_n$, with $K_n$ denoting the Gauss curvature of $\bp_n$. It is known \cite{Hel} that the right-hand side is an element of the Hardy space $\mathcal{H}^1$ so soon as the second fundamental form lies in $L^2$. It follows that $\la_n$ is bounded in $W^{1,(2,1)}$ \cite{BR1}. The analogous $\la_n\in W^{2,(2,1)}$ in dimension 4 must be imposed separately from the onset. \\

In dimension 4, the control of $\la$ is granted by the Paneitz equation. The Paneitz operator in dimension 4 is a fourth-order, conformally invariant differential operator that plays a central role in geometric analysis and conformal geometry. Introduced by Stephen Paneitz \cite{paneitz2008quartic}, it acts on functions and generalises the Laplacian in a conformally invariant way. For a 4-dimensional Riemannian manifold with metric $g$, the Paneitz operator is defined by
$$
P_g\la = -\Delta_g^2\la + d^\star\left( \left( \frac{2}{3} R_g g - 2 \mathrm{Ric}_g \right) \res d\la \right),
$$
Under a conformal change of metric ${g} = e^{2\la}g_0$, the operator transforms as
$$
P_{g}\la = e^{-4\la} P_{g_0}\la,
$$
showing its conformal covariance. It is intimately tied to Q-curvature via the Paneitz equation:
$$
P_{g_0} \la +  Q_{g_0} =  e^{4\la}Q_{{g}}\:,
$$
which is the 4-dimensional equivalent of the Liouville equation in dimension 2. We remind the reader that the Q-curvature $Q_{g_0}$ is explicitly given by\footnote{We refer the reader to \cite{WhatisQ} for an introduction to Q-curvature.}
$$
Q_{g_0}\;=\;-\dfrac{1}{6}\Delta_{g_0}R_{g_0}-\dfrac{1}{2}|\text{Ric}_{g_0}|^2+\dfrac{1}{6}R^2_{g_0}\:.
$$
and similarly for $Q_g$. \\

\noindent
As far as I know, there is no better way to write Q-curvature and this is insufficient to obtain regularity estimates under the critical hypothesis that $\bp\in W^{3,2}\cap W^{1,\infty}$. In that case, $Q_g$ lies in $W^{-2,2}\oplus L^1$, which is too weak to proceed. This is a great difference between the Gauss curvature in 2d and its analogue the Q-curvature in dimension 4: there seems to be no ``good" formulation of $Q_g$. For an interesting discussion on this topic, see \cite{Eastwood}.\\

\noindent
Once imposed that $\la_n\in W^{2,(2,1)}$, it follows by the Sobolev-Lorentz embedding theorem that $\la_n$ is uniformly bounded in $L^\infty$. This guarantees no collapse/explosion along the sequence of immersions.

\setcounter{equation}{0} 
\reset

\section{Proofs of Main Theorems}

\subsection{Regularity: proof of Theorem~\ref{Threg}}

As the Euler-Lagrange equation is quite unmanageable when written as in (\ref{horrib}), the first thing to do is to express it in divergence-form. This can be accomplished by evaluating the Noether field corresponding to translation invariance. More precisely, in \cite{Ber2}, we show 
\begin{Th}\label{cons}
A critical point of $\mathcal{E}$ satisfies the conservation law
$$
d^{\star_g}\bV\;=\;\vec{\mathcal{W}}\;=\;\vec{0}\:,
$$
where the 1-vector-valued 1-form $\bV$ satisfies
\bes
\bV^j\;=\;G^{ij}\nabla_i\bp-\pro\nabla_i\bF^{ij}+\bC^j\:.
\ees
In this expression, we have set
\bes
\left\{\begin{array}{lcl}
G^{ij}&=&\dfrac{1}{2}\bh^{ij}\cdot\Delta_\perp\bH+ |\pro\nabla\bH|^2g^{ij}-2\nabla^j\bH\cdot\pro\nabla^i\bH+\mathcal{O}(|\bh|^4)\\[1ex]
\bF^{ij}&=&-\dfrac{1}{2}\Delta_\perp\bH\,g^{ij}+\mathcal{O}(|\bh|^3) \\[1ex]
\bC^j&=&-2(\bh^{ij}\cdot\nabla_i\bH)\bH+2(\bH\cdot\bh^{ij})\pro\nabla_i\bH\:.\end{array}\right.
\ees
Equivalently, 
\be\label{defVVV}
\bV\;=\;\dfrac{1}{2}d\Delta_\perp\bH+\vec{V}_0\:,
\ee
with
$$
\vec{V}_0^j\;=\;\nabla_k\bA^{kj}+\mathcal{O}\big(|\bh|^2|\pro d\bH|+|\pro d\bH|^2+|\bh|^4\big)\:,
$$
and
$$
\bA^{kj}\;:=\;\mathcal{O}\big(|\bh||\pro d\bH|+|\bh|^3\big)\:.
$$
\end{Th}
$\hfill\square$\\

From the point of view of regularity, which will concern us in this section, a few general observations are first in order. Let $\Omega$ be an open ball in $\R^4$. By hypothesis, we have $\bp\in W^{3,2}\cap W^{1,\infty}(\Omega)$. The induced metric $g_{ij}=\partial_i\bp\cdot\partial_j\bp$ lies in $W^{2,2}\cap L^\infty$. This space is easily seen to be an algebra under pointwise multiplication, which we will use repeatedly. 
Let $B$ be an open proper sub-ball of $\Omega$. Note that the Christoffel symbols satisfy $\Gamma_g\in W^{1,2}(B)\subset L^4(B)$. Since $\partial_{ij}\bp$ lies in $W^{1,2}(B)$, we find the second fundamental form
$$
\bh_{ij}\;=\;\nabla_{ij}\bp\;=\;\partial_{ij}\bp+\mathcal{O}(|\Gamma_g||d\bp|)\:.
$$
The second term on the right-hand side lies in the space
$$
W^{1,2}\cdot (W^{2,2}\cap L^\infty)\;\subset\;W^{1,2}\:.
$$
This shows that $\bh_{ij}$ lies in $W^{1,2}(B)$. In fact, raising indices using the $g$ metric does not affect this fact. In this setting, the energy
$$
\mathcal{E}\;=\;\int\big(|d\bH|_g^2+\mathcal{O}(|\bh|_g^4)\big)d\text{vol}_g
$$
is defined. \\

Given a tensor $T$ with components in $L^{4/3}(B)$, it will be useful for us to note that 
$$
(\nabla_g)T\;=\;\partial T+\mathcal{O}(|\Gamma_g||T|)\;\in\;W^{-1,4/3}\oplus L^4\cdot L^{4/3}\;\subset\;W^{-1,4/3}\oplus L^1\;\subset\;W^{-2,2}\oplus L^1\:,
$$
where we have used the dual Sobolev injection in the last step. The space $W^{-2,2}\oplus L^1$ is easily seen to be preserved under multiplication by a coefficient in $W^{2,2}\cap L^\infty$. Therefore, when $T$ is a tensor in $L^{4/3}(B)$, we can unequivocally write: $\nabla_gT\in W^{-2,2}\oplus L^1(B)$. Similarly if $T$ has components in $L^2(B)$, then 
$$
(\nabla_g)T\;\in\;W^{-1,2}\oplus L^4\cdot L^{2}\;\subset\;W^{-1,2}\oplus L^{4/3}\;\subset\;W^{-1,2}\:.
$$
and this space is preserved under multiplication by a coefficient in $W^{2,2}\cap L^\infty$. Again, we can unequivocally write $\nabla_gT\in W^{-1,2}(B)$. \\

We will impose the smallness hypothesis
\be\label{epsreg}
\Vert\bh\Vert_{L^{4,2}(\Omega)}\;<\;\eps_0\:,
\ee
for some $\eps_0$, which will be adjusted to satisfy our needs.\\
Our setup guarantees that
$$
\Delta_g\bH\;=\;|g|^{-1/2}\partial_i(|g|^{1/2}g^{ij}\partial_j\bH)\;\in\;W^{-1,2}(B)\:.
$$
This implies
$$
\Delta_\perp\bH\;=\;\Delta_g\bH+g^{ai}g^{bj}(\bH\cdot\bh_{ij})\bh_{ab}-2(\nabla_g)_i\big(|\bH|^2g^{ij}-\bH\cdot\bh^{ij}\big)\nabla_j\bp\;\in\;W^{-1,2}\oplus L^{4/3}\;\subset\;W^{-1,2}(B)\:.
$$
More precisely, we have the estimate
\be\label{sikko1}
\Vert \Delta_\perp\bH\Vert_{W^{-1,2}(B)}\;\lesssim\;\Vert \pro d\bH\Vert_{L^2(B)}+\Vert\bh\Vert_{L^{4}(B)}\:,
\ee
where here and throughout the paper, the symbol $\lesssim$ indicates the presence of an irrelevant multiplicative constant involving only the given data on $\Omega$, namely $\eps_0$, the elliptic constants of $g$, and the $W^{2,2}\cap L^\infty$ norm of $g$. \\
Using the above discussion along with (\ref{sikko1}), we use (\ref{defVVV}) to come to the estimate
\be\label{sikko5}
\Vert\vec{V}\Vert_{W^{-2,(2,\infty)}(B)}\;\leq\;\Vert\vec{V}\Vert_{L^1\oplus W^{-2,2}(B)}\;\lesssim\;\Vert\bh\Vert_{L^4(B)}+\Vert\pro d\bH\Vert_{L^2(B)}\:,
\ee
where the first inequality follows trivially from the dual Sobolev injection. \\

A critical point $\bp$ of $\mathcal{E}$ satisfies $d^{\star_g}\bV=\vec{0}$, where $\star_g$ is the Hodge operator relative to the induced metric $g$. We deduce the existence of a 2-form $\bL_0$ satisfying 
\be\label{siko1}
d^{\star_{g}}\bL_0\;=\;\bV\:.
\ee
With (\ref{sikko5}), the latter may be put into Lemma~\ref{lemregL2} from the Appendix to conclude that we can choose $\bL_0$ to be closed and with estimate
\be\label{ppsg1}
\Vert\bL_0\Vert_{W^{-1,(2,\infty)}(B)}\;\lesssim\;\Vert\bh\Vert_{L^4(B)}+\Vert\pro d\bH\Vert_{L^2(B)}\:.
\ee
Now that we have this estimate under our belt, we need to inspect the structure of the vector-field $\bV$ more closely. For reasons that will become clear (and are discussed in detail in \cite{Ber1}), we change $\bL_0$ as follows:
$$
\bL\;:=\;\bL_0-\dfrac{1}{3}d|\bH|^2\stackrel{\wedge}{\wedge}d\bp\:.
$$
This 2-form $\bL$ is also closed, since $\bL_0$ is. Clearly, (\ref{ppsg1}) gives
\be\label{psg1}
\Vert\bL\Vert_{W^{-1,(2,\infty)}(B)}\;\lesssim\;\Vert\bh\Vert_{L^4(B)}+\Vert\pro d\bH\Vert_{L^2(B)}\:.
\ee
Note that
\begin{eqnarray*}
d^{\star_g}\big(d|\bH|^2\stackrel{\wedge}{\wedge}d\bp\big)&=&\nabla_i\big(\nabla^i|\bH|^2\nabla^j\bp-\nabla^j|\bH|^2\nabla^i\bp\big)\\
&=&\Delta_g|\bH|^2\,\nabla^j\bp-\nabla^{ij}|\bH|^2\nabla_i\bp+\bh^{ij}\nabla_i|\bH|^2-4\bH\nabla^j|\bH|^2\:.
\end{eqnarray*}
Hence
\be\label{ikko1}
\dfrac{1}{3}d^{\star_g}\big(d|\bH|^2\stackrel{\wedge}{\wedge}d\bp\big)\stackrel{\!\!\!\cdot}{\res}d\bp\;=\;\Delta_g|\bH|^2\:,
\ee
Moreover, using Codazzi, 
\begin{eqnarray}\label{ikko2}
 \dfrac{1}{3}d^{\star_g}\big(d|\bH|^2\stackrel{\wedge}{\wedge}d\bp\big)\stackrel{\!\!\!\wedge}{\res}d\bp&=&\dfrac{1}{3}\nabla_i|\bH|^2\bh^{ij}\wedge\nabla_j\bp-\dfrac{4}{3}\nabla^j|\bH|^2\bH\wedge\nabla_j\bp\nonumber\\[1ex]
 &=&\dfrac{1}{3}\nabla_i\big((|\bH|^2(\bh^{ij}-4\bH g^{ij})\wedge\nabla_j\bp\big)-\dfrac{1}{3}|\bH|^2\pi_T\nabla_i(\bh^{ij}-4\bH g^{ij})\wedge\nabla_j\bp\nonumber\\[1ex]
 &=&\dfrac{1}{3}\nabla_i\big((|\bH|^2(\bh^{ij}-4\bH g^{ij})\wedge\nabla_j\bp\big)\:.
 \end{eqnarray}
It is shown in \cite{Ber2} that
%
%
\begin{Prop}\label{prop1}
The following identities hold:
$$\bV\stackrel{\!\!\!\cdot}{\res}_gd\bp\;=\;d^{\star_g}d|\bH|^2\qquad\text{and}\qquad\bV\stackrel{\!\!\!\wedge}{\res}_gd\bp\;=\;d^{\star_g}\bT\:,
$$
with 
\begin{eqnarray}\label{sikko4}
\bT_i&=&\dfrac{1}{2}\big(\Delta_\perp\bH+\langle\bH\cdot\bh,\bh\rangle-7|\bH|^2\bH\big)\wedge\nabla_i\bp+2\bH\wedge\pro\nabla_i\bH+2g^{jk}(\bH\cdot\bh_{ij})\bH\wedge\nabla_k\bp+\mathcal{O}(|\bh|^3)\nonumber\\
&=&\mathcal{O}\big(|\Delta_\perp\bH|+|\bH||\pro d\bH|+|\bh|^3\big)\:.
\end{eqnarray}
\end{Prop}
$\hfill\square$\\

\noindent
Using this Proposition along with (\ref{siko1}) and (\ref{ikko1}) gives
\begin{eqnarray*}
d^{\star_{g}}\big(\bL\stackrel{\!\!\!\cdot}{\res}_{g}d\bp\big)&=&d^{\star_{g}}\bL\stackrel{\!\!\!\cdot}{\res}_{g}d\bp\;\;=\;\;\bV\stackrel{\!\!\!\cdot}{\res}_{g}d\bp-\Delta_g|\bH|^2\\
&=&0\:.
\end{eqnarray*}
We deduce the existence of a 2-form $S$ satisfying
\be\label{defS}
d^{\star_{g}} S\;=\;\bL\stackrel{\!\!\!\cdot}{\res}_{g}d\bp\:.
\ee
For reasons that will become clear in the sequel, we ask that
\be\label{deffS}
dS\;=\;2\bL\stackrel{\cdot}{\wedge}d\bp\:,
\ee
which is possible since $\bL$ is closed. \\
From (\ref{psg1}) and the trivial fact that $W^{-1,2}\subset W^{-1,(2,\infty)}$, we have
$$
\Vert dS\Vert_{W^{-1,(2,\infty)}(B)}+\Vert d^{\star_{g}}S\Vert_{W^{-1,(2,\infty)}(B)}\;\lesssim\;\Vert\bh\Vert_{L^4(B)}+\Vert\pro d\bH\Vert_{L^2(B)}\:.
$$
As proved in Lemma~\ref{lemregSR}, it is possible to choose $S$ with the estimate
\be\label{psg6}
\Vert S\Vert_{L^{2,\infty}(B)}\;\lesssim\;\Vert\bh\Vert_{L^4(B)}+\Vert\pro d\bH\Vert_{L^2(B)}\:.
\ee

\medskip

We can of course repeat the same routine with the ambient wedge product in place of the dot product using this time (\ref{ikko2}), namely:
\begin{eqnarray*}
d^{\star_{g}}\big(\bL\stackrel{\!\!\!\wedge}{\res}_{g}d\bp\big)&=&d^{\star_{g}}\bL\stackrel{\!\!\!\wedge}{\res}_{g}d\bp\;\;=\;\;\bV\stackrel{\!\!\!\wedge}{\res}_{g}d\bp-\dfrac{1}{3}\nabla_i\big((|\bH|^2(\bh^{ij}-4\bH g^{ij})\wedge\nabla_j\bp\big)\\
&=&d^{\star_g}\bU\:,
\end{eqnarray*}
where
\be\label{ikko3}
\bU^i\;:=\;\bT^i-\dfrac{1}{3}|\bH|^2(\bh^{ij}-4\bH g^{ij})\wedge\nabla_j\bp\:.
\ee
We deduce the existence of a 2-vector-valued 2-form $\bR$ satisfying
\be\label{defR}
d^{\star_{g}}\bR\;=\;\bL\stackrel{\!\!\!\wedge}{\res}_{g}d\bp-\bU\:.
\ee
For future purposes, we set
\be\label{deffR}
d\bR\;=\;2\bL\stackrel{\wedge}{\wedge}d\bp\:.
\ee
Proposition~\ref{prop1} and the discussion prior (specifically (\ref{sikko1})) provide the following estimate\footnote{we tacitly use again the dual Sobolev injection $L^{4/3}\subset W^{-1,2}$.}
$$
\Vert \bU\Vert_{W^{-1,2}(B)}\;\lesssim\;\Vert\bh\Vert_{L^4(B)}+\Vert\pro d\bH\Vert_{L^2(B)}\:.
$$
This, along with (\ref{psg1}) brought into (\ref{defR})-(\ref{deffR}) yields
$$
\Vert d\bR\Vert_{W^{-1,(2,\infty)}(B)}+\Vert d^{\star_{g}}\bR\Vert_{W^{-1,(2,\infty)}(B)}\;\lesssim\;\Vert\bh\Vert_{L^4(B)}+\Vert\pro d\bH\Vert_{L^2(B)}\:.
$$
We again appeal to Lemma~\ref{lemregSR} to produce a solution $\bR$ with the estimate
\be\label{psg7}
\Vert \bR\Vert_{L^{2,\infty}(B)}\;\lesssim\;\Vert\bh\Vert_{L^4(B)}+\Vert\pro d\bH\Vert_{L^2(B)}\:.
\ee

\bigskip

The next step is far less trivial and is settled in details in \cite{Ber2}. It shows that the potentials $S$ and $\bR$ satisfy a system of equations with an interesting structure:
\begin{Prop}\label{ThSR}
If there exists a closed 2-form $\bL$-system such that the 2-forms $S$ and $\bR$ satisfy
$$
d^\star S=\bL\stackrel{\!\!\!\cdot}{\res}d\bp\quad,\quad dS=2\bL\stackrel{\cdot}{\wedge}d\bp\quad,\quad d^{\star}\bR=\bL\stackrel{\!\!\!\wedge}{\res}d\bp-\bU\quad,\quad dR=2\bL\stackrel{\wedge}{\wedge}d\bp\:,
$$
with $\bU$ as in (\ref{defU}), then the following system holds:
\begin{equation}\label{sysSR20}
\left\{\begin{array}{rcl}
-3d^{\star_g}\bR&=&\vec{\eta}\stackrel{\!\!\!\bullet}{\res}d^{\star_g}\bR+d\bR\stackrel{\!\!\!\bullet}{\res}\vec{\eta}+\vec{\eta}\res d^{\star_g} S-dS\res\vec{\eta}+\mathcal{O}\big(|\bh|^3+|\bh||\pro d\bH|\big)\\
3d^{\star_g} S&=&\vec{\eta}\stackrel{\!\!\!\cdot}{\res}d^{\star_g}\bR-d\bR\stackrel{\!\!\!\cdot}{\res}\vec{\eta}+\mathcal{O}\big(|\bh|^3+|\bh||\pro d\bH|\big)\:,
\end{array}\right. 
\end{equation}
where $\vec{\eta}:=\dfrac{1}{2}d\bp\stackrel{\wedge}{\wedge}d\bp$. 
\end{Prop}
$\hfill\square$\\[2ex]
\noindent
Of interest, of course, is the fact that $\bL$ has been eliminated from the equations. But just as importantly, the ``bad" contributor of $\bU$, namely $\frac{1}{2}\Delta_g\bH\wedge d\bp$, has also disappeared. These are two fundamental advantages of the $(S,\bR)$-system given in Proposition~\ref{ThSR}. Yet another advantage is the possibility to further simplify the system, as hinted by the following result, established in \cite{Ber2}. 
\begin{Prop}
It holds
\bes
d^{\star_g}\big(3\bR+\vec{\eta}\stackrel{\bullet}{\odot}\bR+\vec{\eta}\odot S\big)+d\big(\langle\vec{\eta}\stackrel{\bullet}{,}\bR\rangle+\langle\vec{\eta},S\rangle\big)\;=\;\mathcal{O}\big(|\bh|^3+|\bh||\pro d\bH|+|\bh||\bR|+|\bh||S|\big)\:.
\ees
\end{Prop}
$\hfill\square$

\bigskip

The regularity proof truly begins now. We first impose a smallness hypothesis:
\be\label{epsreg1}
\Vert \bh\Vert_{L^{4,2}(\Omega)}\;<\;\eps_0\:.
\ee
We have the estimates
\begin{eqnarray}\label{psg8}
\big\Vert|\bh|(|\bR|+|S|)\big\Vert_{L^{4/3,2}(B)}&\lesssim&\Vert\bh\Vert_{L^{4,2}(B)}\big(\Vert\bR\Vert_{L^{2,\infty}(B)}+\Vert S\Vert_{L^{2,\infty}(B)}\big)\nonumber\\
&\lesssim&\eps_0\big(\Vert\bh\Vert_{L^4(B)}+\Vert\pro d\bH\Vert_{L^2(B)}\big)\:,
\end{eqnarray}
where in the last step we have used (\ref{psg6}) and (\ref{psg7}). \\[2ex]
\noindent
For notational convenience, let $\bq:=\langle\vec{\eta}\stackrel{\bullet}{,}\bR\rangle+\langle\vec{\eta},S\rangle$ and let $\bv$ denote the 1-form on the right-hand side of the first identity in the above system:
\be\label{psg9}
d^{\star_g}\big(3\bR+\vec{\eta}\stackrel{\bullet}{\odot}\bR+\vec{\eta}\odot S\big)+d\bq\;=\;\bv\;=\;\mathcal{O}\big(|\bh|^3+|\bh||\pro d\bH|+|d^\star\vec{\eta}||\bR|+|d^\star\vec{\eta}||S|\big)\:.
\ee
Per (\ref{psg8}) and (\ref{epsreg1}), it comes
\be\label{psg10}
\Vert\bv\Vert_{L^{4/3,2}(B)}\;\lesssim\;\eps_0\big(\Vert\bh\Vert_{L^4(B)}+\Vert\pro d\bH\Vert_{L^2(B)}\big)\:.
\ee
Applying the codifferential to both sides of (\ref{psg9}) yields
$$
\Delta_g\bq\;=\;d^{\star_g}\bv\:\qquad\text{on}\:\:\Om\:,
$$
or equivalently 
\be\label{defq}
\partial_i\big(|g|^{1/2}g^{ij}\partial_j\bq\big)\;=\;\partial_i\big(|g|^{1/2}g^{ij}\bv_j\big)\:,
\ee
The coefficients $|g|^{1/2}g^{ij}$ are uniformly elliptic, non-degenerate, and belong to the space $W^{2,2}\cap L^\infty\subset \text{VMO}\cap L^\infty$. Interior elliptic estimates\footnote{linearly interpolate to Lorentz norms.} are provided e.g. in \cite{DiFazio}:
\be\label{estimq0}
\Vert d\bq\Vert_{L^{4/3,2}(\Om_k)}\;\lesssim\;\Vert \bv\Vert_{L^{4/3,2}(\Om)}+k^a\Vert\bq\Vert_{L^{2,\infty}(B)}\:,
\ee
with $B_k$ has the same centre as $B$ and its radius rescaled by a factor $k\in(0,1)$; and $a>0$ is some constant. 
By definition, $\bq=\mathcal{O}(|S|+|\bR|)$, so that by (\ref{psg6}) and (\ref{psg7}), we find
\be\label{vac1}
\Vert \bq\Vert_{L^{2,\infty}(B)}\;\lesssim\;\Vert\bh\Vert_{L^4(B)}+\Vert\pro d\bH\Vert_{L^2(B)}\:.
\ee
Injected into (\ref{estimq0}) and using (\ref{psg10}), we come to 
\be\label{estimq}
\Vert d\bq\Vert_{L^{4/3,2}(\Om_k)}\;\lesssim\;(\eps_0+k^a)\big(\Vert\bh\Vert_{L^4(B)}+\Vert\pro d\bH\Vert_{L^2(B)}\big)\:.
\ee

\bigskip

In order to proceed with the proof, we need to show how to control geometry by $\vec{q}$. This was done in \cite{Ber2}:
\begin{Prop}
It holds
\begin{eqnarray*}
2\Delta_\perp\bH&=&\pro d^{\star_g}(\bq\res d\bp)+\mathcal{O}\big(|\bh|^3+|\bh||\pro d\bH|+|\bh||\bR|+|\bh||S|\big)\:.
\end{eqnarray*}
\end{Prop}
\hfill$\square$\\

Bringing into this proposition (\ref{vac1}), (\ref{estimq}), as well as (\ref{psg6}) and (\ref{psg7}), proves that
\bes
\Vert\Delta_\perp\bH\Vert_{L^{4/3,2}(\Om_k)}\;\lesssim\;(\eps_0+k^a)\big(\Vert\bh\Vert_{L^4(\Om)}+\Vert\pro d\bH\Vert_{L^2(\Om)}\big)\:.
\ees
Since 
$$
\Delta_g \bH\;=\;\Delta_\perp\bH-(\bH\cdot\bh^{ij})\bh_{ij}+2\nabla_i\big(|\bH|^2g^{ij}-\bH\cdot\bh^{ij}\big)\nabla_j\bp\:,
$$
and $L^{4/3,2}\subset W^{-1,2}$, we may convert the previous estimate into 
\be\label{estimdelH}
\Vert\Delta_g\bH\Vert_{W^{-1,2}(\Om_k)}\;\lesssim\;(\eps_0+k^a)\big(\Vert\bh\Vert_{L^4(\Om)}+\Vert\pro d\bH\Vert_{L^2(\Om)}\big)\:.
\ee
As in the Appendix, we may again appeal to \cite{Mcintosh} (for our bounded metric) to deduce the estimate
\be\label{estimdH}
\Vert\pro d\bH\Vert_{L^2(\Om_k)}\;\lesssim\;(\eps_0+k^b)\big(\Vert\bh\Vert_{L^4(\Om)}+\Vert\pro d\bH\Vert_{L^2(\Om)}\big)\:,
\ee
for some $b>0$. \\

%
%
%
To proceed further, we have to find a suitable bound for $\Vert\bh\Vert_{L^4(B_k)}$. This may be achieved as follows. Introduce the 2-vector-valued 1-form
$$
\vec{\frak{h}}^a\;:=\;\bh^{ai}\wedge\nabla_i\bp\:.
$$
Using the Codazzi identity, it is elementary to check that
$$
d^{\star_g}\vec{\frak{h}}\;=\;\nabla_a\bh^{ai}\wedge\nabla_i\bp\;=\;4\pro\nabla^i\bH\wedge\nabla_i\bp
$$
and
\begin{eqnarray*}
d\vec{\frak{h}}&=&\nabla^a(\bh^{bi}\wedge\nabla_i\bp)-\nabla^b(\bh^{ai}\wedge\nabla_i\bp)\\
&=&\pi_T(\nabla^a\bh^{bi}-\nabla^b\bh^{ai})\nabla_i\bp+2\bh^{bi}\wedge\bh_i^a\\
&=&R^{abij}\nabla_i\bp\wedge\nabla_j\bp+2\bh^{bi}\wedge\bh^a_i\\
&=&\mathcal{O}(|\bh|^2)\:.
\end{eqnarray*}
From these computations, we find that
$$
\Vert d\vec{\frak{h}}\Vert_{L^2(\Om_k)}+\Vert d^{\star_{g}}\vec{\frak{h}}\Vert_{L^2(\Om_k)}\;\lesssim\;\Vert\bh\Vert^2_{L^4(\Om_k)}+\Vert\pro d\bH\Vert_{L^2(\Om_k)}\;\lesssim\;\eps_0\Vert\bh\Vert_{L^4(\Om_k)}+\Vert\pro d\bH\Vert_{L^2(\Om_k)}\:,
$$
where we have again used the smallness hypothesis (\ref{epsreg}). Although our metric $g$ is not continuous, it is bounded, and Gaffney type inequalities are still at reach for the exponent $p=2$ (see \cite{Mcintosh}). In particular, we have now
\begin{eqnarray*}
\Vert D\vec{\frak{h}}\Vert_{L^2(\Om_k)}&\lesssim&\Vert d\vec{\frak{h}}\Vert_{L^2(\Om_k)}+\Vert d^{\star_{g}}\vec{\frak{h}}\Vert_{L^2(\Om_k)}+\Vert \vec{\frak{h}}\Vert_{L^2(\Om_k)}\\
&\lesssim&(\eps_0+k^c)\Vert\bh\Vert_{L^4(\Om_k)}+\Vert\pro d\bH\Vert_{L^2(\Om_k)}\:,
\end{eqnarray*}
for some $c>0$. 
Using the Sobolev embedding theorem, we deduce an estimate for $\vec{\frak{h}}$ in $L^4$, and correspondingly, we arrive at 
\be\label{estimdh0}
\Vert d\bh\Vert_{L^2(\Om_{k^2})}+\Vert \bh\Vert_{L^4(\Om_{k^2})}\;\lesssim\;(k^b+\eps_0)\Vert\bh\Vert_{L^4(\Om_k)}+\Vert\pro d\bH\Vert_{L^2(\Om_k)}\:.
\ee
Combining this with (\ref{estimdH}) shows that
\bes
\Vert \bh\Vert_{L^4(\Om_{k^2})}+\Vert\pro d\bH\Vert_{L^2(\Om_{k^2})}\;\lesssim\;(\eps_0+k^b)\big(\Vert\bh\Vert_{L^4(\Om)}+\Vert\pro d\bH\Vert_{L^2(\Om)}\big)\:.
\ees
From standard growth estimates (cf. \cite{GM}), we reach the following Morrey-type decay:
\be\label{morreydec}
\Vert \bh\Vert_{L^4(\Om_{r})}+\Vert\pro d\bH\Vert_{L^2(\Om_{r})}\;\lesssim\;r^\beta\:,
\ee
for some $\beta>0$. In turn, once imported into (\ref{estimdelH}), the latter yields
\bes
\Vert\Delta_\perp\bH\Vert_{L^{4/3,2}(\Om_{r})}\;\lesssim\;r^\beta\:,
\ees
Since 
$$
\Delta_g \bH\;=\;\Delta_\perp\bH-(\bH\cdot\bh^{ij})\bh_{ij}+2\nabla_i\big(|\bH|^2g^{ij}-\bH\cdot\bh^{ij}\big)\nabla_j\bp\:,
$$
the previous two estimates give now
$$
\big\Vert\partial_j(a^{ij}\partial_i\bH)\big\Vert_{L^{4/3,2}(B_r)}\;\lesssim\;r^\beta
$$
with coefficients $a^{ij}=|g|^{1/2}g^{ij}$ that are uniformly elliptic and lie in the space $W^{2,2}\cap L^{\infty}\subset W^{1,4}$. The latter combined with (\ref{morreydec}) yields in particular
\begin{eqnarray*}
\big\Vert a^{ij}\partial_{ij}\bH\big\Vert_{L^{4/3,2}(B_r)}&\lesssim&r^\beta+\Vert Da\Vert_{L^4(B_r)}\Vert d\bH\Vert_{L^2(B_r)}\;\;\lesssim\;\;r^\beta\:.
\end{eqnarray*}
We may now call upon standard $L^p$-theory for second-order uniformly elliptic equations in non-divergence form with coefficients in VMO, such as the results in \cite{Chiarenza}. Interpolating to Lorentz norms, we obtain
\be\label{morrey0}
\Vert\bH\Vert_{W^{2,(4/3,2)}(B_{r})}\;\lesssim\;r^\beta\:,
\ee
and thus in particular
\be\label{morrey}
\Vert\Delta_{\text{flat}}\bH\Vert_{L^{4/3,2}(B_{r})}\;\lesssim\;r^\beta\:.
\ee

\medskip

In the next step, we introduce the maximal function
$$\mathcal{M}_{3-\beta}f\;:= \;\sup_{1>r>0}\,r^{-1-\beta}\Vert{f}\Vert_{L^1(B_r)}.$$

\noindent
By Jensen's inequality, we have 
$$\mathcal{M}_{3-\beta}f\;\lesssim\; \sup_{1>r>0}\, r^{-\beta}\Vert f\Vert _{L^{4/3,2}(B_r)}.$$
Using the Morrey decay (\ref{morrey}), we obtain
\be\label{cecon2}
\big\Vert\mathcal{M}_{3-\beta}\Delta_{\text{flat}}\bH\big\Vert_{L^\infty(B_r)} \;\lesssim\;1\;\quad\quad\forall\,\, r<1.  
\ee

\noindent
For a locally integrable function $f$ on $\mathbb{R}^4$, the {\it Riesz potential $\mathcal{I}_1$} of $f$ is defined by the convolution
$$(\mathcal{I}_1 *f)(x):=\int_{\mathbb{R}^4} f(y)|x-y|^{-1} dy\:.$$

 \noindent
We will now use the following interesting result from \cite{Ada}.  
\begin{Prop}  \label{nbd}
If $0<\lambda\leq 4$, $1<p<\lambda$, and $f\in L^p(\mathbb{R}^4)$ with $\mathcal{M}_{\lambda/p} f\in L^\infty(\Omega)$, $\Omega\subset \mathbb{R}^4$, then
$$
\Vert\mathcal{I}_1f\Vert_{L^s(\Omega)} \lesssim \Vert \mathcal{M}_{\lambda/p }f\Vert_{L^\infty(\Omega)}^{ p/\lambda} \Vert{f}\Vert_{L^p(\Omega)}^{1- p/\lambda} 
$$
where $1/s=1/p-1/\lambda$.
\end{Prop}
$\hfill\square$
\vskip3mm
\noindent
Putting $p\in(1,4/3]$, and $\lambda=(3-\beta)p$ in Proposition \ref{nbd} and using (\ref{cecon2}), we find 
\be\label{cecon}
\big\Vert\mathcal{I}_1\Delta_{\text{flat}}\bH\big\Vert_{L^{s}(B_r)}\; \lesssim \; \big\Vert\Delta_{\text{flat}}\bH\big\Vert_{L^{p}(B_r)}^{(2-\beta)/(3-\beta)}\:,
\ee
where
$$
s:= p\left( \frac{3-\beta}{2-\beta} \right)\;>\;2\:.
$$
Choosing $p<4/3$ but very close to $4/3$, we arrive at
\bes
\big\Vert\mathcal{I}_1\Delta_{\text{flat}}\bH\big\Vert_{L^{s}(B_r)}\; \lesssim \; \big\Vert\Delta_{\text{flat}}\bH\big\Vert_{L^{4/3,2}(B_r)}^{(2-\beta)/(3-\beta)}\:,
\ees
for some $s>2$.
Introducing (\ref{morrey}) yields that $\big\Vert\mathcal{I}_1\Delta_{\text{flat}}\bH\big\Vert_{L^s(B_r)}$ is bounded, which finally gives
$$
\Vert d \bH\Vert_{L^s(B_r)}\; \lesssim \;1\:.
$$

\bigskip

The regularity of $\bH$ has thus improved to $W^{1,s}$ for some $s>2$, and so $\bH\in L^{s^*}$ with $s^*>4$ denoting the conjugate Sobolev exponent of $s$. In order to update this information back to the immersion, we might be tempted to use the fact that $\Delta_g\bp=4\bH$. Applying derivative to both sides yields $\Delta_{\text{Hodge}_g}d\bp=4d\bH$. Since the right-hand side belongs to $L^s$ for some $s>2$, it follows from \cite{Mcintosh} that $d\bp$ lies in $W^{2,s}$ (it suffices to ensure that $2<s<p^H$, where $p^H$ is given in \cite{Mcintosh}). In other words, $\bp\in C^{0,\alpha}$. This regularises the metric $g$ and the entire problem may now be revisited with $g$ being H\"older continuous. A standard bootstrapping procedure eventually confirms that $\bp$ is a smooth immersion, and geometric estimates follow accordingly. This completes the proof of Theorem~\ref{Threg}.

\subsection{Palais-Smale sequences: proof of Theorem~\ref{ThPS}}

We have seen that critical points of our energy satisfy 
$$
d^{\star_g}\bV\;=\;\vec{0}\:,
$$
for some vector-field $\bV$ which we have discussed at length. For a (non-degenerate) metric of the form $g=e^{2\la}g_0$ we have seen in the paragraph preceding (\ref{siko1}) that
$$
d^{\star_{g_0}}(e^{2\la}\bV)\;=\;\vec{0}\:.
$$
In this section we consider a sequence of immersions $\bp_n\in W^{3,2}\cap W^{1,\infty}(B)$ such that the corresponding sequence of induced metrics $g_n$ remains within the same conformal class: $g_n=e^{2\la_n}g_0$. We demand that $(\la_n)$ be uniformly bounded in $W^{2,(2,1)}$. As far as the fixed metric $g_0$ is concerned, we suppose that $g_0\in C^\infty$ is uniformly elliptic and non-degenerate.
We further impose the following {\it Palais-Smale}-type condition\footnote{Prime indicates the dual space.}:
\be\label{defPS}
d^{\star_{g_0}}(e^{2\la_n}\bV_n)\:\longrightarrow\:\vec{0}\qquad\text{strongly in}\:\:\:(W^{3,2}\cap W^{1,\infty})'\:,
\ee
The vector-field $\bV_n$ is associated to $\bp_n$ via Theorem~\ref{cons}. \\
%

From the Sobolev-Lorentz injection, we know that $W^{3,(2,1)}$ embeds continuously into $W^{3,2}\cap W^{1,\infty}$.
 In particular, (\ref{defPS}) gives
\be\label{defPS2}
d^{\star_{g_0}}(e^{2\la_n}\bV_n)\:\in\:(W^{3,(2,1)})'\;=\;W^{-3,(2,\infty)}
\ee
with uniform estimates. 
From the distributional characterisation of dual Sobolev spaces (interpolated to Sobolev-Lorentz norms) given for example in Theorem 3.10 from \cite{AF}, and the smoothness of $g_0$, there exists a 1-form $\vec{y}_n\in W^{-2,(2,\infty)}$ such that 
$$
\partial_j\vec{y}_n^j\;=\;\partial_j\big(|g_0|^{1/2}g_0^{ij}e^{2\la_n}(\bV_n)_i\big)\qquad\text{and}\qquad \Vert\vec{y}_n\Vert_{W^{-2,(2,\infty)}}\;\leq\;\Big\Vert \partial_j\big(|g_0|^{1/2}g_0^{ij}e^{2\la_n}(\bV_n)_i\big)\Big\Vert_{W^{-3,(2,\infty)}}\:.
$$
It follows in particular that $\vec{y}_n$ converges to $\vec{0}$ strongly in $W^{-2,(2,\infty)}$. Define next $\vec{Y}_n$ via the identity
$$
\vec{y}_n^j\;:=\;|g_0|^{1/2}g_0^{ij}e^{2\la_n}(\vec{Y}_n)_i\:.
$$
The following continuous injection is elementary to verify
$$
W^{2,(2,1)}\cdot W^{-2,(2,\infty)}\;\subset\;W^{-2,(2,\infty)}\:. 
$$
Since  $\la_n$ is assumed to be uniformly bounded in $W^{2,(2,1)}$, we find that $\vec{Y}_n$ converges to $\vec{0}$ strongly in $W^{-2,(2,\infty)}$. 
In addition, we have by construction
$$
d^{\star_{g_0}}(e^{2\la_n}\vec{Y}_n)\;=\;d^{\star_{g_0}}(e^{2\la_n}\bV_n)\:.
$$
We deduce the existence of a 2-form $\bL_n$ such that
\be\label{cor01}
d^{\star_{g_0}}\bL_n\;=\;e^{2\la_n}(\vec{V}_n-\vec{Y}_n)
\ee
We have seen that $e^{2\la_n}\vec{Y}_n$ lies in $W^{-2,(2,\infty)}$ by construction and converges to $0$ strongly in $W^{-2,(2,\infty)}$. In addition, from the way $\bV_n$ is constructed from $\bp_n$, it comes easily that $\bV_n$ lies in the space $W^{-2,(2,\infty)}$, with uniform estimate (\ref{sikko5}). As $e^{4\la_n}$ lies in $W^{2,(2,1)}$ which continuously preserves $W^{-2,(2,\infty)}$ under multiplication, we deduce from (\ref{cor01}) that
\be\label{cor1}
d^{\star_{g_0}}\bL_n\;\in\;W^{-2,(2,\infty)}\qquad\text{uniformly in $n$}\:.
\ee
Lemma~\ref{lemregL2} guarantees that we can choose $\bL_n$ to be closed and moreover the following estimate holds
\be
\Vert \bL_n\Vert_{W^{-1,(2,\infty)}}\;\lesssim\;\Vert\vec{Y}_n\Vert_{W^{-2,(2,\infty)}}+\Vert\vec{V}_n\Vert_{W^{-2,(2,\infty)}}\:.
\ee
 As the right-hand side is uniformly bounded, we can extract a subsequence, still denoted $(\bL_n)_n$ which converges weak* in $W^{-1,(2,\infty)}$ to some limit $\bL$:
 \be\label{cure}
 \bL_n\; \overset{*}{\rightharpoonup} \bL\qquad\text{in}\quad W^{-1,(2,\infty)}\:.
 \ee
On the other hand because $\bp_n$ is uniformly bounded in $W^{3,2}$, we can extract a subsequence converging weakly in $W^{3,2}$ to some limit $\bp$. The Rellich-Kondrachov theorem enables us to further extract a subsequence converging strongly to $\bp$ in $\bigcap_{p<4}W^{2,p}$. In particular, $d\bp_n$ converges strongly to $d\bp$ in $\bigcap_{p<4}W^{1,p}$. This space is easily seen to be an algebra, and it follows that $g_n$ converges to $g$ strongly in $\bigcap_{p<4}W^{1,p}$, where $g$ denotes the pull-back metric of the limit immersion $\bp$. We deduce that $\la_n$ converges strongly in $\bigcap_{p<4}W^{1,p}$ to some $\la$ with the property that $g=e^{2\la}g_0$. As $\la_n$ is uniformly bounded in $W^{2,(2,1)}\subset L^\infty$, the limit $\la$ is bounded. That insures $\bp$ is an immersion.  
Since $W^{-1,(2,\infty)}$ and $W^{1,(2,1)}$ are in functional duality, it follows that\footnote{As several metrics appear in this section, we append a subscript to the $\res$ symbol, to remind the reader that the first-order contraction operator requires raising indices, and thus depends on the metric.}
\bes
\left\{\begin{array}{lcl}
\bL_n\stackrel{\!\!\!\cdot}{\res}_{g_0} d\bp_n&\longrightarrow&\bL\stackrel{\!\!\!\cdot}{\res}_{g_0} d\bp\\[1ex]
\bL_n\stackrel{\!\!\!\wedge}{\res}_{g_0} d\bp_n&\longrightarrow&\bL\stackrel{\!\!\!\wedge}{\res}_{g_0} d\bp  
\end{array}\right.\qquad\text{in}\:\:\:\mathcal{D}'\:.
\ees
We obtain
\bes
\left\{\begin{array}{lclcl}
d^{\star_{g_0}}\big(\bL_n\stackrel{\!\!\!\cdot}{\res}_{g_0} d\bp_n\big)&\longrightarrow&d^{\star_{g_0}}\big(\bL\stackrel{\!\!\!\cdot}{\res}_{g_0} d\bp\big)&=&d^{\star_{g_0}}\bL\stackrel{\!\!\!\cdot}{\res}_{g_0} d\bp\\[1ex]
d^{\star_{g_0}}\big(\bL_n\stackrel{\!\!\!\wedge}{\res}_{g_0} d\bp_n\big)&\longrightarrow&d^{\star_{g_0}}\big(\bL\stackrel{\!\!\!\wedge}{\res}_{g_0} d\bp\big) &=&d^{\star_{g_0}}\bL\stackrel{\!\!\!\wedge}{\res}_{g_0} d\bp  
\end{array}\right.\qquad\text{in}\:\:\:\mathcal{D}'\:.
\ees
Because $g=e^{2\la}g_0$, it is elementary to check that
$$
d^{\star_{g_0}}\big(\bL_n\stackrel{\!\!\!\cdot}{\res}_{g_0} d\bp_n\big)\;=\;d^{\star_{g_0}}\bL_n\stackrel{\!\!\!\cdot}{\res}_{g_0}d\bp_n\;=\;e^{4\la_n}d^{\star_{g_n}}\bL_n\stackrel{\!\!\!\cdot}{\res}_{g_n}d\bp_n\:,
$$
and similarly
$$
d^{\star_{g_0}}\big(\bL_n\stackrel{\!\!\!\wedge}{\res}_{g_0} d\bp_n\big)\;=\;e^{4\la_n}d^{\star_{g_n}}\bL_n\stackrel{\!\!\!\wedge}{\res}_{g_n}d\bp_n\:.
$$
Owing to the identities following Proposition~\ref{prop1}, the latter may further be recast as
\be\label{rusdem}
\left\{\begin{array}{lcl}
d^{\star_{g_0}}\big(e^{2\la_n}d|\bH_n|^2\big)&\longrightarrow&d^{\star_{g_0}}\bL\stackrel{\!\!\!\cdot}{\res}_{g_0} d\bp\\[1ex]
d^{\star_{g_0}}\big(e^{2\la_n}\bU_n\big)&\longrightarrow&d^{\star_{g_0}}\bL\stackrel{\!\!\!\wedge}{\res}_{g_0} d\bp  
\end{array}\right.\qquad\text{in}\:\:\:\mathcal{D}'\:,
\ee
where the 1-form $\bU_n$ is constructed from the immersion $\bp_n$ as per Proposition~\ref{prop1}:
\be\label{psg13}
(\bU_n)_j\;:=\;\dfrac{1}{2}\Delta_{g_n}\bH_n\wedge\nabla_j\bp_n+\mathcal{O}(|\bh_n|^3+|\bh_n||\pro d\bH_n|\big)\:.
\ee

As $\bh_n$ is uniformly bounded in $W^{1,2}$, we can further extract a subsequence of immersions with corresponding sequence of second fundamental form converging weakly in $W^{1,2}$ and strongly in $\bigcap_{p<4}L^p$ to the second fundamental form of the limit immersion. In particular, since $\bh_n$ and $d\bH_n$ are uniformly bounded in $L^4$ and $L^2$ respectively, while $\bh$ converges strongly to $\bh$ in $\bigcap_{p<4}L^p$ and $d\bH_n$ converges weakly to $d\bH$ in $L^2$, it follows that $d|\bH_n|^2$ converges to $d|\bH|^2$ weakly in $\bigcap_{s<4/3}L^s$. As $e^{2\la_n}$ is uniformly bounded in $W^{2,2}\cap L^\infty$ and converges strongly to $e^{2\la}$ in $\bigcap_{p<4}W^{1,p}$, we deduce that $e^{2\la_n}d|\bH_n|^2$ converges weakly to $e^{2\la}|d\bH|^2$ in $\bigcap_{s<4/3}L^s$. The first item in (\ref{rusdem}) then show that
$$
d^{\star_{g_0}}(e^{2\la}d|\bH|^2)\;=\;d^{\star_{g_0}}\bL\stackrel{\!\!\!\cdot}{\res}_{g_0}d\bp\qquad\text{in}\:\:\:\mathcal{D}'\:.
$$
This is of course equivalent to 
\be\label{PSeq1}
d^{\star_{g}}d|\bH|^2\;=\;d^{\star_{g}}\bL\stackrel{\!\!\!\cdot}{\res}_{g}d\bp\qquad\text{in}\:\:\:\mathcal{D}'\:.
\ee

\medskip
\noindent
By the same token, the remainder term in (\ref{psg13}) converges weakly to its limit in $\bigcap_{s<4/3}L^s$. Furthermore, because $|g_n|^{1/2}g_n^{ij}$ converges strongly to $|g|^{1/2}g^{ij}$ in $W^{1,p}$ for $p<4$, it follows that $|g_n|^{1/2}g_n^{ij}\partial_j\bH_n$ converges weakly in $\bigcap_{q<2} L^q$ to $|g|^{1/2}g^{ij}\partial_j\bH$. In turn, this shows $|g_n|^{1/2}\Delta_{g_n}\bH_n$ converges to $|g|^{1/2}\Delta\bH$ weakly in $\bigcap_{q<2}W^{-1,q}$. We have thus on one hand the sequence $|g_n|^{-1/2}$ which is uniformly bounded in $W^{2,2}\cap L^\infty$ and converges strongly to $|g|^{-1/2}$ in $\bigcap_{p<4}W^{1,p}$, and on the other hand a sequence which is uniformly bounded in $W^{-1,2}$ and which converges weakly in $\bigcap_{q<2}W^{-1,q}$. The product sequence will then converge to limit product weakly in the space $\bigcap_{q<2}W^{-1,q}$. This follows from the injections
$$
\bigcap_{p<4}W^{1,p}\cdot \bigcap_{q'>2}W^{1,q'}\;\subset\;\bigcap_{q'>2}W^{1,q'}\;=\;\bigcap_{q<2}(W^{-1,q})'\:,
$$
and
$$
W^{-1,2}\cdot\bigcap_{q'>2}W^{1,q'}\;\subset\;\bigcap_{p'>4/3}W^{-1,p'}\;=\;\bigcap_{p<4}(W^{1,p})'\:.
$$
Accordingly, $\Delta_{g_n}\bH_n$ converges weakly to $\Delta_g\bH$ in $\bigcap_{q<2}W^{-1,q}$. In an analogous fashion, it now follows from (\ref{psg13}) that $e^{2\la_n}\bU_n$ converges weakly to $e^{2\la}\bU$ in $\bigcap_{q<2}W^{-1,q}$, where $\bU$ is the construct corresponding to the limit immersion $\bp$. Bringing this new information in the second item from (\ref{rusdem}), we find that
$$
d^{\star_{g_0}}\big(e^{2\la}\bU\big)\;=\;d^{\star_{g_0}}\bL\stackrel{\!\!\!\wedge}{\res}_{g_0} d\bp\qquad\text{in}\:\:\:\mathcal{D}'\:.
$$
Equivalently,
\be\label{PSeq2}
d^{\star_{g}}\bU\;=\;d^{\star_{g}}\bL\stackrel{\!\!\!\wedge}{\res}_{g} d\bp\qquad\text{in}\:\:\:\mathcal{D}'\:.
\ee

The system formed by (\ref{PSeq1}) and (\ref{PSeq2}) was studied extensively in \cite{Ber2}, where it is known as the ``$\bL$-system". In that paper it is shown that when the $\bL$-system holds, there exists a divergence-free and traceless symmetric 2-tensor $T$ such that
$$
d^{\star_g}\bV\;=\;\langle T,\bh\rangle\:.
$$
This completes the proof of Theorem~\ref{ThPS}.

\bigskip

\setcounter{equation}{0} 
\reset

\renewcommand{\theequation}{A.\arabic{equation}}
\renewcommand{\theTh}{A.\arabic{Th}}
\renewcommand{\theProp}{A.\arabic{Prop}}
\renewcommand{\theLm}{A.\arabic{Lm}}
\renewcommand{\theCo}{A.\arabic{Co}}
\renewcommand{\theRm}{A.\arabic{Rm}}
\renewcommand{\theequation}{A.\arabic{equation}}
\setcounter{equation}{0} 
\reset
\appendix
\section{Appendix}

\subsection{Notational Conventions}\label{nota}
\reset

We fix a basis $\{e^i\}_{i=1,\cdots,N}$ of $\mathbb{R}^N$ and an orientation. Suppose $\mathbb{R}^N$ is equipped with a metric $\frak{g}$. Let $\Lambda^k(\R^N)$ be the corresponding vector space of $k$-forms. We will need component expansions of certain operations involving forms. A typical element $A\in\Lambda^p$ may be expanded in tensorial notation as
$$
A\;=\;\dfrac{1}{p!}\sum_{\text{all indices}}A_{i_1\ldots i_p}\,e^{i_1}\wedge\ldots\wedge e^{i_p}\:,
$$
where $A_{i_1\ldots i_p}$ is totally antisymmetric. For notational convenience, we will render this as
$$
A\;\equiv\;A_{i_1\ldots i_p}\:.
$$
Forms are tensors: we raise and lower indices as usual using the metric $\frak{g}$. Of course, as forms are anti-symmetric, the relative position of indices is essential: we follow standard conventions. \\
The inner product on forms will be denoted
$$
\langle A,B\rangle_\frak{g}\;=\;\dfrac{1}{p!}A^{i_1\ldots i_p}B_{i_1\ldots i_p}\qquad\text{for}\quad A\in\Lambda^p\:\:,\:\: B\in\Lambda^p\:.
$$
The bracket symbol $[\ldots]$ will denote \underline{scaled} total permutations. We conventionally set the wedge product to be
$$
A\wedge B\;\equiv\;\dfrac{(p+q)!}{p!q!}A_{[i_1\ldots i_p}B_{j_1\ldots j_q]}\qquad\text{for}\quad A\in\Lambda^p\:\:,\:\: B\in\Lambda^q\:.
$$
The {\it interior multiplication} (restricted to our interest) is the bilinear operation $\res\;:\;\Lambda ^p(\mathbb{R}^N)\times \Lambda^{p-1}(\mathbb{R}^N)\rightarrow \Lambda^{1}(\mathbb{R}^N)$ which satisfies in particular
$$
A\res B\;\equiv\;\dfrac{1}{(p-1)!}A^{i_1\ldots i_p}B_{i_1\ldots i_{p-1}}\qquad\text{for}\quad A\in\Lambda^p\:\:,\:\: B\in\Lambda^{p-1}\:.
$$
\noindent
The {\it first-order contraction operator}   $\bullet: \Lambda^2(\mathbb{R}^N)\times \Lambda^q(\mathbb{R}^N)\rightarrow \Lambda^{q}(\mathbb{R}^N)$ is defined as follows. For a $2$-form ${A}$ and a $1$-form ${B}$, we set 
$${A}\bullet {B}:={A}\res B\:.$$
For a $p$-form ${B}$ and a $q$-form ${C}$, we have
$${A}\bullet ({B}\wedge{C}):=({A}\bullet {B})\wedge {C}+ (-1)^{pq}({A}\bullet{C})\wedge{B}\:.$$
This operation will be used for forms with values in $\Lambda^p(\R^{m\ge5})$.\\
We will also need an operation taking a pair a two-forms to a two-form, namely:
$$
A\odot B\;\equiv\;A^{i_2j}B_{j}^{\:\:i_1}-A^{i_1j}B_{j}^{\:\:i_2}\qquad\text{for}\quad A, B\in\Lambda^2\:.
$$
\noindent
Being given an orientation and a metric, we can equip our manifold with the volume form $\epsilon\in\Lambda^N$ via $\epsilon^{i_1\ldots i_N}:=|\frak{g}|^{-1/2}\mbox{sign} \begin{pmatrix}
1 & 2 & \cdots & N\\
i_1 &i_2 & \cdots & i_N
\end{pmatrix}$. With it, we define the Hodge star operator $\star_\frak{g}:\Lambda^p\rightarrow\Lambda^{N-p}$ as usual. In components, we have
$$
\star_\frak{g}A\;\equiv\;\dfrac{1}{p!}\,\epsilon^{i_1\ldots i_{N-p}j_1\ldots j_p}A_{j_1\ldots j_p}\qquad\text{for}\quad A\in\Lambda^p\:.
$$
Note that
$$
\star_\frak{g}\star_\frak{g} A=(-1)^{p(N-p)}A\:.
$$

\noindent
If a connection $\nabla$ compatible with the metric $\frak{g}$ is provided, we define in accordance with the Leibniz rule, the covariant exterior derivative $d$ of a $p$-form $A$ in components as
$$dA\;:=\;(p+1)\nabla_{[i_0}A_{i_1...i_p]}\;\in\;\Lambda^{p+1}\:.$$
We also define the codifferential $d^{\star_\frak{g}}:=-\star_\frak{g} d\star_\frak{g}$, whose components take the expression
$$
d^{\star_\frak{g}} A=\nabla^{j}A_{ji_1...i_{p-1}}\;\in\;\Lambda^{p-1}\:.
$$
For any $p$-form $A$, it holds
$$
d^2A\;=\;0\;=\;(d^{\star_\frak{g}})^2A\:.
$$

\medskip

\noindent
We will work within the context of multivector-valued forms. Our objects will be elements of $\Lambda^p(\R^4,\Lambda^q(\R^{m\ge5}))$. To distinguish operations taking place in the parameter space $(\R^4,g)$ from operations taking place in the ambient space $(\R^{m},g_{\text{Eucl}})$, we will superimpose the operators defined above with the convention that the one on the bottom acts on parameter space, while the one at the top acts on ambient space. 
Typically, consider the two-vector-valued two-form $\vec{\eta}\in\Lambda^2(\R^4,\Lambda^2(\R^m))$ whose components are given by
$$
\vec{\eta}_{ij}\;=\;\nabla_i\bp\wedge\nabla_j\bp\:.
$$
We will write more succinctly
$$
\vec{\eta}\;=:\;\dfrac{1}{2}d\bp\stackrel{\wedge}{\wedge}d\bp\:.
$$
Similar quantities are defined following the same logic. Naturally, for multivector 0-forms, the bottom symbol will be omitted, while for 0-vector-valued forms, the symbol on the top will be omitted instead. The notation $\langle\cdot\,,\cdot\rangle$ is reserved for $\langle\cdot\,,\cdot\rangle_g$.

\subsection{Sobolev-Lorentz spaces}

For the reader's convenience, we recall in this section the fundamentals of Lorentz spaces. Detailed accounts may be found in \cite{BL} and \cite{Ta1}.\\[-1ex]

For a real-valued measurable function $f$ on an open subset of $U\subset \mathbb{R}^n$, its belonging to a Lorentz space is determined by a condition involving the non-decreasing rearrangement of $|f|$ on the interval $(0,|U|)$, where $|U|$ denotes the Lebesgue measure of $U$. The non-increasing rearrangement $f^*$ of $|f|$ is the unique positive continuous function from $(0,|U|)$ into $\mathbb{R}$ which is non-increasing and satisfies
\bes
\Big|\big\{x\in U\;\big|\;|f(x)|\ge s\big\}\Big|\:=\:\Big|\big\{t\in(0,|U|)\;\big|\;f^*(t)\ge s\big\}\Big|\:.
\ees
\noindent
If $p\in(1,\infty)$ and $q\in[1,\infty]$, the Lorentz space $L^{p,q}(U)$ is the set of measurable functions $f:U\rightarrow\mathbb{R}$ for which
\bes
\int_{0}^{\infty}\,\big(t^{1/p}f^*(t)\big)^q\,\dfrac{dt}{t}\:<\:\infty\qquad\text{if}\:\:\:q<\infty\:,
\ees 
or
\bes
\sup_{t>0}\,t^{1/p}f^*(t)\:<\:\infty\qquad\text{if}\:\:\:q=\infty\:.
\ees \\
A complete norm on $L^{p,q}(U)$ is given by
\bes
\Vert f\Vert_{L^{p,q}(U)}\:=\:\big\Vert \,t^{1/p}f^{**}\big\Vert_{L^q([0,\infty),\,dt/t)}\qquad\text{where}\qquad f^{**}(t)\,:=\,\dfrac{1}{t}\,\int_{0}^{t}\,f^*(\tau)\,d\tau\:.
\ees
One verifies that
\bes
L^p\;=\;L^{p,p}\:,
\ees
and that $L^{p,\infty}$ is the weak-$L^p$ Marcinkiewicz space. \\[1ex]
Moreover, we have the inclusions
\bes
L^{p,1}\;\subset\;L^{p,q'}\;\subset\;L^{p,q''}\;\subset\;L^{p,\infty}\qquad\text{for}\qquad 1<q'<q''<\infty\:;
\ees
and if $U$ has finite measure, there holds for all $q$ and $q'$
\bes
L^{p',q'}(U)\;\subset\;L^{p,q}(U)\qquad\text{whenever}\qquad p<p'\:.
\ees
More precisely,
\be\label{Lorinc}
\Vert f\Vert_{L^{p,q}(U)}\;\;\lesssim\;\;|U|^{\frac{p'-p}{p'p}}\,\Vert f\Vert_{L^{p',q'}(U)}\:.
\ee

\medskip
\noindent
Finally, if $q<\infty$, the space $\,L^{\frac{p}{p-1},\frac{q}{q-1}}$\, is the dual of $L^{p,q}$.\\[1.5ex]
Similarly to Lebesgue spaces, Lorentz spaces obey a pointwise multiplication rule and a convolution product rule. More precisely, for $\,1<p_1,p_2<\infty\,$ and $\,1\le q_1,q_2\le\infty$, there holds
\bes
L^{p_1,q_1}\,\times\,L^{p_2,q_2}\:\subset\:L^{p,q}\qquad\text{with}\qquad\left\{\begin{array}{lcl}
p^{-1}&=&p_1^{-1}+\,p_2^{-1}\\[1ex]
q^{-1}&=&q_1^{-1}+\,q_2^{-1}\:,
\end{array}\right.
\ees
and
\bes
L^{p_1,q_1}\,*\,L^{p_2,q_2}\:\subset\:L^{p,q}\qquad\text{with}\qquad\left\{\begin{array}{lcl}
p^{-1}&=&p_1^{-1}+\,p_2^{-1}-\,1\\[1ex]
q^{-1}&=&q_1^{-1}+\,q_2^{-1}\:.
\end{array}\right.
\ees\\
\noindent
An interesting feature of Lorentz spaces is that they are interpolation spaces between Lebesgue spaces. \\

\noindent
Lorentz spaces offer the possibility to sharpen the classical Sobolev embedding theorem. More precisely, it can be shown (see \cite{BL} and \cite{Ta3}) that
\bes
W^{k,q}(\mathbb{R}^m)\:\subset\:L^{p,q}(\R^m)
\ees
is a continuous embedding so long as
\bes
1\;\le\;q\;\le\;p\;<\;\infty\qquad\text{and}\qquad \dfrac{k}{m}\;=\;\dfrac{1}{q}\,-\,\dfrac{1}{p}\:.
\ees

\medskip
The Sobolev-Lorentz spaces are defined analogously to the ``standard" Sobolev spaces, but with the Lebesgue norms replaced by Lorentz norms. The  Sobolev-Lorentz space $\,W^{m,(p,q)}(\di)\,$ consists of all locally summable functions $u$ on $\di$ such that $\,D^\al u\,$ exists in the weak sense and belongs to the Lorentz space $\,L^{p,q}(\di)\,$ for each multiindex $\,\al\,$ with $\,|\al|\leq m$. The norm
\bes
\Vert u\Vert_{W^{m,(p,q)}}\;:=\;\sum_{|\al|\leq m}\Vert D^\al u\Vert_{L^{p,q}}
\ees
clearly makes $\,W^{m,(p,q)}(\di)\,$ into a Banach space. \\
The Sobolev-Lorentz embedding theorem \cite{Ta3} may be written as the continuous injection
$$
W^{1,(p,q)}(\R^m)\;\subset\;L^{p^*,q}(\R^m)\:,
$$
where $p^*$ is the conjugate Sobolev exponent of $p$, while $q\in[1,\infty]$. Of exception is the remarkable embedding
$$
W^{1,(m,1)}(\R^m)\;\subset\;L^\infty(\R^m)\:.
$$

\subsection{Hodge-type estimates}

Let $\Om\subset\mathbb{R}^4$ be a ball with smooth boundary. We will suppose that the induced metric $g_{ij}$ is uniformly elliptic, non-degenerate, and with bounded coefficients. In this situation, it is possible to obtain Hodge decompositions with ``natural" estimates for a suitable range of parameters. This is explained in details in \cite{Mcintosh}. In particular, any $k$-form $A\in\Lambda^k(\Om)$, for $k\in\{1,2,3\}$ can be decomposed on $\Om$ {\it à la} Hodge:
\bes
A\;=\;da+d^{\star_{g}} b+\om\:,
\ees
where $\om$ is a harmonic $k$-form: $d\om=0=d^{\star_{g}}\om$, $a$ is a $(k-1)$-form, and $b$ is a $(k+1)$-form. The tangential part of $a$ and the normal part of $b$ can be chosen to vanish on the boundary $\partial\Om$, which we write $a_T=0=b_N$, following the notation and definitions given in \cite{Iwa}\footnote{We deliberately elude the technical aspects of defining traces of forms on the boundary $\partial B$, and refer the reader to \cite{Mcintosh}. See also \cite{Iwa} and \cite{Li}}. When the entries of $A$ are in $L^p(\Om)$, for some $p\in (p_H,p^H)$, then the following control is guaranteed:
\be\label{eq1}
\Vert a\Vert_{W^{1,p}(\Om)}+\Vert b\Vert_{W^{1,p}(\Om)}\;\lesssim\;\Vert A\Vert_{L^p(\Om)}\:.
\ee
The exponents $p_H$ and $p^H$ are known to satisfy $p_H<2<p^H$. In particular, the above estimate can be interpolated to the Lorentz space $L^{(2,1)}$. Similarly, Calderon-Zygmund type $L^p$ estimates for the second-order problem $\Delta_{\text{Hodge}_g}A\in L^p$ -- with appropriate boundary conditions -- are available in the same range of parameter $p$. See again \cite{Mcintosh}. \\

\begin{Lm}\label{lemregL2}
Let $g$ be uniformly elliptic, non-degenerate, and bounded. Suppose that a 2-form $L$ satisfies on $B$:
$$
d^{\star_{g}}L\;=\;F\;\in\;W^{-2,(2,\infty)}(B)\:.
$$
We can choose $L$ such that 
$$
\Vert L\Vert_{W^{-1,(2,\infty)}(B)}\;\lesssim\;\Vert F\Vert_{W^{-2,(2,\infty)}(B)}\:.
$$
\end{Lm}
{\bf Proof.} Fix $L$ to be closed. Affecting $L$ by a harmonic form if necessary, we may without loss of generality suppose that $L=dv$, for some 1-form $v$ with boundary condition $v_T=0$. Let $\varphi\in W^{1,(2,1)}$ be a test 2-form. We can decompose $\varphi=d\psi+b$ with $d^{\star_{g}}b=0$ and the boundary condition $\psi_T=0$. Integrating by parts, we have
\begin{eqnarray*}
\int_B\langle L,\varphi\rangle_{g}d\text{vol}_{g}&=&\int_B\langle dv,\varphi\rangle_{g}d\text{vol}_{g}\;\;=\;\;\int_B\langle v,d^{\star_{g}}\varphi\rangle_{g}d\text{vol}_{g}\\
&=&\int_B\langle v,d^{\star_{g}}d\psi\rangle_{g}d\text{vol}_{g}\;\;=\;\;\int_B\langle d^{\star_{g}}dv,\psi\rangle_{g}d\text{vol}_{g}\\
&=&\int_B\langle F,\psi\rangle_{g} d\text{vol}_{g}\:.
\end{eqnarray*}
Because $W^{-1,(2,\infty)}$ is in duality with $W^{1,(2,1)}$, this shows the proof will be done if we can establish that 
\be\label{estpsi}
\Vert\psi\Vert_{W^{2,(2,1)}(B)}\;\leq\;\Vert\varphi\Vert_{W^{1,(2,1)}(B)}\:.
\ee
As $d\psi$ appears in the Hodge decomposition of $\varphi$, we are free to demand that $d^{\star_{g}}\psi=0$. In this case, $\psi$ satisfies the equation
$$
\Delta_{\text{Hodge}_{g}}\psi\;=\;d^{\star_{g}}\varphi\;\in\;L^{2,1}(B)\qquad\text{with boundary condition}\quad\psi_T=0\:.
$$
Per our previous discussion and the aforementioned results from \cite{Mcintosh}, (\ref{estpsi}) follows and the proof is complete.

$\hfill\blacksquare$\\


\begin{Rm} In the case when $F$ belongs to the smaller space $W^{2,2}\oplus L^1$, it might be possible to improve the statement of Lemma~\ref{lemregL2} by proceeding as follows. In the proof, suppose instead the weaker condition $\varphi\in W^{1,2}$. As we did, we can estimate $\varphi$ in $W^{2,2}$. Generalised Bourgain-Brezis-type results \cite{VS} and \cite{Mazya} show that any form in the critical Sobolev space $W^{2,2}$ is bounded up to a closed form:
$$
\psi\;=\;\psi+d\omega\:\qquad\text{with}\qquad\Vert\psi\Vert_{W^{2,2}\cap L^\infty}\;\leq\;\Vert\psi\Vert_{W^{2,2}}\:.
$$
If this result can be localised to the ball $B$ with boundary condition $\omega_T=0$, the following identities arise
\begin{eqnarray*}
\int_B\langle L,\varphi\rangle_{g}d\text{vol}_{g}&=&\int_B\langle F,\psi\rangle_{g} d\text{vol}_{g}\\
&=&\int_B\langle F,\psi_0+d\omega\rangle_{g} d\text{vol}_{g}\;\;=\;\;\int_B\langle F,\psi_0\rangle_{g} d\text{vol}_{g}\:,
\end{eqnarray*}
where we have used that $F$ is co-closed. This line of argument would improve the estimate in Lemma~\ref{lemregL2} to
$$
\Vert L\Vert_{W^{-1,2}(B)}\;\lesssim\;\Vert F\Vert_{W^{-2,2}\oplus L^1(B)}\:.
$$
Once imported into the regularity proof for critical points of $\mathcal{E}$, it would be possible to replace the smallness hypothesis $\Vert\bh\Vert_{L^{4,2}(\Om)}<\eps_0$ by the somewhat more natural $\Vert\bh\Vert_{L^{4}(\Om)}<\eps_0$.
\end{Rm}

$\hfill\square$\\

\begin{Lm}\label{lemregSR}
Let $g$ be as in the previous lemma. Let $ S$ be a 2-form with $d S$ and $d^{\star_{g}} S$ bounded in $W^{-1,(2,\infty)}(\Om)$, for $\Om\subset\mathbb{R}^4$ a ball. Then, affecting $ S$ by a harmonic 2-form if necessary, it holds
$$
\Vert  S\Vert_{L^{2,\infty}(\Om)}\;\lesssim\;\Vert d S\Vert_{W^{-1,(2,\infty)}(\Om)}+\Vert d^{\star_{g}} S\Vert_{W^{-1,(2,\infty)}(\Om)}\:.
$$
\end{Lm}
{\bf Proof.} Affecting $ S$ by a harmonic 2-form if necessary, we may without loss of generality suppose 
$$
 S\;=\;du+d^{\star_{g}} v\qquad\text{with}\quad u_T=0=v_N\:.
$$
Let $A$ be a test 2-form. We introduce the Hodge decomposition $A=da+d^{\star_{g}}b+\mu$, for some harmonic 2-form $\mu$ and $a_T=0=b_N$ on the boundary of $B$. Suppose that $A\in L^{2,1}(B)$ so that by the Sobolev-Lorentz inequality and (\ref{eq1})
\be\label{dac1}
\Vert a\Vert_{W^{1,(2,1)}(\Om)}+\Vert b\Vert_{W^{1,(2,1)}(\Om)}\;\lesssim\;\Vert A\Vert_{L^{2,1}(\Om)}\:.
\ee
Next, we integrate by parts
\begin{eqnarray*}
\int_B\langle  S,A\rangle_{g}d\text{vol}_{g}&=&\int_B\langle du+d^{\star_{g}}v,A\rangle_{g}d\text{vol}_{g}\\
&=&\int_B\langle u,d^{\star_{g}}A\rangle_{g}d\text{vol}_{g}+\int_B\langle v,dA\rangle_{g}d\text{vol}_{g}\\
&=&\int_B\langle u,d^{\star_{g}}da\rangle_{g}d\text{vol}_{g}+\int_B\langle v,dd^{\star_{g}}b\rangle_{g}d\text{vol}_{g}\\
&=&\int_B\langle d^{\star_{g}}du,a\rangle_{g}d\text{vol}_{g}+\int_B\langle dd^{\star_{g}}v,b\rangle_{g}d\text{vol}_{g}\\
&=&\int_B\langle d^{\star_{g}} S,a\rangle_{g}d\text{vol}_{g}+\int_B\langle d S,b\rangle_{g}d\text{vol}_{g}\:.
\end{eqnarray*}
Since $W^{1,(2,1)}$ and $W^{-1,(2,\infty)}$ are in formal duality, (\ref{dac1}) put into the latter confirms that $ S$ lies in $L^{2,\infty}$ with estimate
$$
\Vert S\Vert_{L^{2,\infty}(B)}\;\leq\;\Vert d S\Vert_{W^{-1,(2,\infty)}(\Om)}+\Vert d^{\star_{g_0}} S\Vert_{W^{-1,(2,\infty)}(\Om)}\:,
$$
as desired.

$\hfill\blacksquare$\\

\bibliographystyle{plain}
\bibliography{references}

\end{document}